%% file: main.tex
\documentclass[oneside]{amsart}

\usepackage{control/packages}
\usepackage{control/docaesthetic}
\usepackage{control/environments}
\usepackage{control/macros}
\usepackage{control/paperdata}

\usepackage{subfiles}

\usepackage[style=numeric]{biblatex}
\begin{document}

\maketitle
\subfile{content/frontmatter/abstract}

\tableofcontents

\subfile{content/frontmatter/introduction}

\subfile{content/mainmatter/prelim}

\subfile{content/mainmatter/complexes}

\subfile{content/mainmatter/complexfunctors}

\subfile{content/mainmatter/complexmonads}

\subfile{content/mainmatter/catrelmorph}

\subfile{content/mainmatter/nerves}

\subfile{content/mainmatter/pointlikes}

\subfile{content/mainmatter/moduli}

\subfile{content/mainmatter/effectiveness}

\subfile{content/mainmatter/transfer}

\subfile{content/mainmatter/reversal}

\subfile{content/mainmatter/contextspec}

\subfile{content/mainmatter/conclusion}

\printbibliography
\hrulefill
\end{document}

%% file: content/frontmatter/abstract.tex
%

\begin{abstract}
We introduce a general unifying framework for the investigation of pointlike sets. 
The pointlike functors are considered as distinguished elements of a certain lattice of subfunctors of the power semigroup functor;
in particular, we exhibit the pointlike functors as the fixed points of a closure operator induced by an antitone Galois connection between this lattice of functors and the lattice of pseudovarieties. 
Notably, this provides a characterization of pointlikes which does not mention relational morphisms. 
Along the way, we formalize various common heuristics and themes in the study of pointlike sets. 
As an application, we provide a general method for transferring lower bounds for pointlikes along a large class of continuous operators on the lattice of pseudovarieties.
\end{abstract}

%

%% file: content/frontmatter/introduction.tex
%

\section{Introduction}

This paper presents an abstract framework for the study of pointlike sets.
This framework---which utilizes category- and lattice-theoretic tools---provides a uniform language for describing, generating conjectures for, and proving results about pointlikes.
The construction of this framework follows the trend of ``turning theorems into definitions'' (see \cite{QTHEORY})---various ubiquitous heuristics and motifs in the study of pointlikes are axiomatized and tied together to form the foundations of a general theory of pointlike sets.

\subsection{Background}
Given a finite semigroup $S$, let $\Po(S)$ denote the semigroup of non-empty subsets of $S$ under the inherited operation given by
    \[
        X \cdot Y = \setst{xy}{x \in X, y \in Y}
    \]
for all non-empty subsets $X$ and $Y$ of $S$, and let $\sing(S)$ denote the subsemigroup of $\Po(S)$ consisting of the singletons.

A morphism $\varphi : S \rightarrow T$ extends to a morphism
    \[
        \ext{\varphi} : \Po(S) \arw \Po(T)
            \quad \text{given by} \quad
        (X)\ext{\varphi} = \setst{(x)\varphi}{x \in X}.
    \]
Equipping the object map $\Po$ with this action on morphisms yields a functor
    \[
        \Po : \FSGP \arw \FSGP
    \]
which creates monomorphisms, regular epimorphisms, and isomorphisms.

For each pseudovariety $\pvV$, there is a \textit{$\pvV$-pointlikes functor}
    \[
        \PV : \FSGP \arw \FSGP
    \]
which is a subfunctor of $\Po$ with the property that a finite semigroup $S$ belongs to $\pvV$ if and only if $\PV(S) = \sing(S)$.\footnotemark
\footnotetext{The notation for $\Po$ is due to it being the pointlikes functor for the trivial pseudovariety $\pvtriv$.}
Pointlike functors also create monomorphisms, regular epimorphisms, and isomorphisms.

For a given $S$, a non-empty set $X \in \Po(S)$ belongs to $\PV(S)$ if for any relational morphism of the form $\rho : S \rlm V$ with $V \in \pvV$ there exists some element $v \in V$ for which $X \subseteq (v)\rho^\inv$.
The constituent sets of $\PV(S)$ are said to be \textit{pointlike} with respect to $\pvV$ (or \textit{$\pvV$-pointlike} for short).

A pseudovariety $\pvV$ is said to have \textit{decidable pointlikes} if there is an algorithm which produces the $\pvV$-pointlike subsets of any finite semigroup given as input.
Decidability of pointlikes implies decidability of membership by way of the aforementioned fixed-point property, but the converse is false: 
pseudovarieties with decidable membership but undecidable pointlikes are given in \cite{UNDECIDABLE, EXTPARPERM}.
Intuitively, whereas membership is boolean, $\pvV$-pointlikes ``measure the essential differences'' between a given semigroup and those which belong to $\pvV$.\footnotemark
\footnotetext{As an illustration, if $G$ is a finite group, then $\PL{\pvA}(G) = \Po(G)$ (where $\pvA$ is the pseudovariety of finite aperiodic semigroups).}

Before we give an overview of the concepts and results of this paper, we briefly review existing results on pointlikes.\footnotemark
\footnotetext{See \cite{STEINBERG-PL-SURVEY} for a dedicated survey.}
The initial motivation for considering pointlikes was the complexity problem,\footnotemark
\footnotetext{See \cite[Chapter~4]{QTHEORY}.}
and so it should come as no surprise that the first pseudovarieties shown to have decidable pointlikes were $\pvA$ (due to the first author in \cite{PL-APERIODIC-KARSTEN}) and $\pvG$ (due to Ash's proof of the Type II conjecture in \cite{ASHTYPEII}, which confirmed a description conjectured by Rhodes and the first author in \cite{PLG-CONJ}).

There have been a number of alternative proofs of the decidability of aperiodic pointlikes (see \cite{PRODEXP, MERGEDECOMP}); moreover, the tools used in the aperiodic case have been generalized to obtain decidability of pointlikes for various pseudovarieties defined by restrictions on subgroups (see \cite{PL-BEYOND, PL-VARDETGRP-BEN}).
The other Green’s-relation-trivial pseudovarieties have decidable pointlikes as well: the cases of $\pvR$ and $\pvL$ were established by Almeida and Silva in \cite{HYP-R-ALM} (see also \cite{PL-RJ-ALM}); and the case of $\pvJ$ was established first by Steinberg in \cite{PL-JOINS-BEN} and later by Almeida \textit{et~al.} in \cite{PL-RJ-ALM}.
In addition to $\pvJ$-pointlikes (and other more general results considered momentarily), in \cite{PL-JOINS-BEN} Steinberg proves a number of pointlike decidability results, including for the pseudovarieties $\pvSL$, $\pvN$, $\pvD$, $\pvK$, and $\mathbb{L}\pvtriv$.\footnotemark
\footnotetext{These pseudovarieties are semilattices, nilpotent semigroups, delay, reverse delay, and locally trivial semigroups, respectively.}

There have also been a number of pseudovarieties of groups which have been shown to have decidable pointlikes.
In \cite{PL-ABELIAN} Delgado shows that the pseudovariety $\pvAb$ of abelian groups has decidable pointlikes; this was followed a few years later by \cite{PL-ABELIAN-ANY}, wherein any decidable pseudovariety of abelian groups is shown to have decidable pointlikes.
Other pseudovarieties of groups which have decidable pointlikes include the pseudovarieties of nilpotent groups \cite{TAME-NILGROUPS} and $p$-groups for any prime $p$ \cite{INEVITABLE-BEN}.

Of particular interest are ``transfer results'' for pointlikes, which are results of the form \textit{``if $\pvV$ has decidable pointlikes, then $(\pvV)\alpha$ has decidable pointlikes''}---possibly contingent on $\pvV$ satisfying some conditions---where $\alpha$ is an operator on the lattice of pseudovarieties.
The framework of this paper was developed with transfer results in mind---in particular, it provides a formal setting in which results of this form may be pursued.

Several transfer results have been established.
The first such results, due to Steinberg in \cite{PL-JOINS-BEN}, involve transfer along the join operation on the lattice of pseudovarieties by way of his ``slice theorem''. 
Notably, he establishes decidability of pointlikes for $\pvW \vee \pvV$ whenever $\pvW$ is locally finite with computable relatively free semigroups and $\pvV$ has decidable pointlikes. 
Additionally, he gives conditions which, when satisfied by a pseudovariety $\pvV$ with decidable pointlikes, imply decidability of pointlikes for $\pvW \vee \pvV$ when $\pvW$ is one of $\pvN$, $\pvD$, $\pvK$, $\mathbb{L}\pvtriv$, $\pvJ$, and $\pvACom$ (aperiodic and commutative semigroups).\footnotemark
\footnotetext{Note that, with some exceptions, each case requires different conditions on $\pvV$. For details, see \cite{PL-JOINS-BEN}.}

Other transfer results concern semidirect products. 
In \cite{PL-DELAYTHM}, Steinberg proves a generalization of Tilson’s derived category theorem (see \cite{DERCAT-TILSON}) which extends pointlikes to the context of finite categories in order to give a characterization of $\pvV \ast \pvW$ pointlikes.
This theorem is then used to show that $\pvV \ast \pvD$ has decidable pointlikes if and only if $\pvV$ does (a generalization of Tilson's delay theorem).

Recent transfer results for semidirect products have come from the theory of regular languages (for reasons discussed below)---for instance, Place and Zeitoun show in \cite{SEPCOVGRP} that $\pvJ \ast \pvH$-pointlikes are decidable whenever $\pvH$ is a pseudovariety of groups with decidable pointlikes.
Additionally, the results of Place \textit{et~al.} in \cite{COVSEPMODPRED} provide conditions on a pseudovariety $\pvV$ with decidable pointlikes whose satisfaction implies the decidability of $\pvV \ast \pvAb$-pointlikes.\footnotemark
\footnotetext{These algebraic translations of the results of \cite{SEPCOVGRP, COVSEPMODPRED} are taken from \cite{STEINBERG-PL-SURVEY}.}

As alluded to regarding the last two results, pointlikes are also relevant to the theory of regular languages.
In the language-theoretic context, pointlikes take the form of the equivalent \textit{covering problem}, which was formulated by Place and Zeitoun in \cite{COVPROB}. 
Also relevant is the \textit{separation problem}, which was shown by Almeida in \cite{ALGPROBPV} to be equivalent to the problem of computing two-element pointlike sets.
Relevant papers on these topics include \cite{SEPGOINGHIGHER, SEPLOCTEST, SEPFOLOG, SEPDOTDEPTH, SEPADDSUCC}.
Despite its relevance, the language-theoretic perspective is beyond the scope of this paper; here we adopt an explicitly algebraic point of view. 
However, it would be interesting to translate the framework developed here into language-theoretic terms.

\subsection{Overview of paper}
Section \ref{section:prelim} establishes notational conventions and briefly reviews various preliminary concepts from category theory, order theory, and finite semigroup theory.
All of the material covered therein (with the possible exception of certain notational and terminological conventions) is standard.

Section \ref{section:complexes} is concerned with an ``object-level'' description of pointlikes---that is, we ask: \textit{given a finite semigroup $S$, what possible values can $\PV(S)$ take for some pseudovariety $\pvV$?}
Our answer is based on the observation that $\PV(S)$ is a subsemigroup of $\Po(S)$ which
\begin{itemize}
    \item contains $\sing(S)$ as a subsemigroup, and which
    \item is closed under taking non-empty subsets of its members; that is, if $X \in \PV(S)$ and $X_0 \subseteq X$ with $X_0 \neq \varnothing$, then $X_0 \in \PV(S)$ as well.
\end{itemize}
Transmuting these properties into axioms yields the notion of a \textit{semigroup complex} (Definition \ref{com:defn:complex}), which is a pair $(S, \cK)$ consisting of a finite semigroup $S$ and a subsemigroup $\cK$ of $\Po(S)$ which contains the singletons and which is closed under taking non-empty subsets of its members.
Their name is due to the fact that these are abstract simplicial complexes whose vertex set is a finite semigroup, and whose faces are also a finite semigroup under the inherited multiplication;\footnotemark for this reason, if $(S, \cK)$ is a semigroup complex then $S$ and $\cK$ are called the \textit{vertex} and \textit{face semigroups}, respectively.
\footnotetext{In categorical terms, semigroup complexes are precisely the semigroup objects in the category of finite abstract simplicial complexes.}

The various face semigroups of semigroup complexes whose vertices are a given semigroup $S$ are called \textit{$S$-complexes}, and the set of $S$-complexes obtains the structure of a complete lattice via a closure operator on the lattice of subsemigroups of $\Po(S)$ (\ref{com:env:scomlattice}).
Moreover, this provides the action on objects of a functor from the category of finite semigroups to the category of complete finite lattices and join-preserving maps, and the category of semigroup complexes is recoverable from this functor via the Grothendieck construction (\ref{com:env:underlyingfunctors}).

Section \ref{section:complexfunctors} builds on the preceding section to characterize the pointlike functors.
To this end, we introduce the notion of a \textit{complex functor} (Definition \ref{cf:defn:complexfunctor}), which is a subfunctor $\cC$ of $\Po$ which preserves regular epimorphisms and which has the property that $(S, \cC(S))$ is a semigroup complex for every finite semigroup $S$.
The collection $\CF$ of complex functors carries the structure of a complete lattice, which is largely inherited pointwise from the lattices of $S$-complexes as $S$ ranges over all finite semigroups (\ref{cf:env:cflattice}).

A finite semigroup $S$ is said to be a \textit{fixed point} of a complex functor $\cC$ if $\cC(S) = \sing(S)$.
The set of fixed points of $\cC \in \CF$ is denoted by $\Fix(\cC)$, i.e.,
    \[
        \Fix(\cC) 
            =
        \setst{S \in \FSGP}{\cC(S) = \sing(S)}.
    \]
The fixed points of a complex functor form a pseudovariety, and the map
    \[
        \Fix : \CF \arw \PVAR,
    \]
where $\PVAR$ denotes the lattice of pseudovarieties, is antitone and takes joins to meets (Proposition \ref{FIRSTPUNCHLINE}).
Consequently, $\Fix$ has an (antitone) upper adjoint---this will turn out in Section \ref{section:pointlikes} to be the map $\PLF$ sending pseudovarieties to their respective pointlike functors (Theorem \ref{MAINTHEOREM}), yielding a Galois connection
    \begin{equation*}
            \begin{tikzcd}
            \CF \arrow[rrr, "\Fix"', bend right=15, two heads, ""{name=B,above}] 
            & & &
            \PVAR^\op. \arrow[lll, "\PLF"', bend right=15, hook', ""{name=A,below}]
            \ar[from=A, to=B, phantom, "\ddash"]
        \end{tikzcd}
    \end{equation*}
This provides a ``relational-morphism-free'' characterization of $\PV$ as the largest complex functor whose pseudovariety of fixed points contains $\pvV$.

Returning to our linear outline, Section \ref{section:complexmonads} is motivated by the well-known fact that each pointlike functor $\PV$ may be naturally equipped with the structure of a monad $(\PV, \sigma_\pvV, \mu_\pvV)$; 
where at each finite semigroup $S$ the component of the unit transformation $\sigma_\pvV : \id{\FSGP} \Rightarrow \PV$ is the singleton embedding
    \[
        \sigma_{\pvV, S} 
            =
        \{-\} : S \longinj \PV(S)
            \quad \text{given by} \quad
        x \longmapsto \{x\}
    \]
and the component of the multiplication $\mu_\pvV : \PV^2 \Rightarrow \PV$ is the union map
    \[
        \mu_{\pvV, S} 
            =
        \bigcup (-) : \PV^2(S) \longsurj \PV(S)
            \quad \text{given by} \quad
        \cX \longmapsto \bigcup_{X \in \cX} X.
    \]

A complex functor which admits a monad structure of this form is called a \textit{complex monad} (Definition \ref{cm:defn:complexmonad}).
The ubiquitous "apply then union and iterate until closed" technique (see \cite{PL-APERIODIC-KARSTEN, PRODEXP, PL-BEYOND, PL-RJ-ALM, PL-VARDETGRP-BEN}) appears here as a closure operator on $\CF$ which sends each $\cC$ to its \textit{monad completion} $\moncom{\cC}$ (Definition \ref{cm:defn:moncom}), which is defined at a finite semigroup $S$ as a sequential colimit over a diagram 
    \begin{equation*}
        \begin{tikzcd}
            \cflevel{\cC}{0}(S)     \ar[r, hook, "i"]
        &
            \cflevel{\cC}{1}(S)     \ar[r, hook, "i"]
        &
            \cdots  \ar[r, hook, "i"]
        &
            \cflevel{\cC}{n}(S)     \ar[r, hook, "i"]
        &
            \cflevel{\cC}{n+1}(S)     \ar[r, hook, "i"]
        &
            \cdots
        \end{tikzcd}
    \end{equation*} 
where $\cflevel{\cC}{n+1}(S)$ is the minimal $S$-complex which contains $\bigcup \cX$ as an element whenever $\cX \in \cC(\cflevel{\cC}{n}(S))$; that is, it is minimal such that the map
    \[
        \bigcup (-) : \cC(\cflevel{\cC}{n}(S)) \arw \cflevel{\cC}{n+1}(S)
    \]
is well-defined.
This equips the collection $\CM$ of complex monads with the structure of a complete lattice.
Moreover, monad completion preserves fixed points.
Altogether, this yields a commutative triangle 
    \begin{equation*}
    \begin{tikzcd}
            \CF
            \ar[rrrr, shift right=1, bend left=15, ""{name=B,above}, two heads, "\moncom{(-)}"']
            \ar[dddrr, shift left=1, bend right=15, ""{name=C}, two heads, "\Fix"]
            &       &       &       &
            \CM
            \ar[llll, shift right=3, bend right=15, ""{name=A,below}, hook', "i"']
            \ar[dddll, shift right=1.5, bend left=15, two heads, ""{name=E}, "\Fix"']
            \ar[from=A, to=B, phantom, "\scriptstyle{\ddash}"]
        \\
            &       &       &       &
        \\
            &       &       &       &
        \\
            &       &
            \PVAR^\op
            \ar[uuull, shift left=3, bend left=15, ""{name=D}, hook, "\PLF"]
            \ar[uuurr, shift right=3, bend right=15, ""{name=F}, hook', "\PLF"']
            &       &
        \ar[from=C, to=D, phantom, "\scriptstyle{\ddash}" rotate=135]
        \ar[from=E, to=F, phantom, "\scriptstyle{\ddash}" rotate=-135]
    \end{tikzcd}
    \end{equation*}
of Galois connections between complete lattices (Proposition \ref{SECONDPUNCHLINE}).\footnotemark
\footnotetext{Although the identity of the upper adjoint to $\Fix$ is not established in Proposition \ref{SECONDPUNCHLINE} (it is not established until Theorem \ref{MAINTHEOREM}), we have previously spoiled that surprise in the introduction, and hence we continue to spoil it here.}

Once semigroup complexes, complex functors, and complex monads each stand under their own weight, we turn to the task of linking this machinery to the problem they are meant to help describe.
This task begins in Section \ref{section:catrelmorph}, which concerns the category whose objects are relational morphisms and whose arrows are the evident pairs of morphisms between domains and codomains. 

This leads into Section \ref{section:nerves}, which defines a crucial functor from this category of relational morphisms to the category of semigroup complexes.
The object map of this functor sends a relational morphism $\rho : S \rlm T$ to its \textit{nerve}, which is the semigroup complex $(S, \Nrv(S , \rho, T))$, where the face semigroup consists of all sets $X \in \Po(S)$ for which there exists some $t \in T$ such that $X \subseteq (t)\rho^\inv$; i.e., such that every member of $X$ is related to $t$ under $\rho$.

Section \ref{section:pointlikes} unites the material of the sections so far and establishes that the map sending pseudovarieties to their pointlike functors is upper adjoint to $\Fix$.

Once this machinery is developed, we pivot towards ``practical'' concerns.
Establishing decidability of $\pvV$-pointlikes for some $\pvV \in \PVAR$ requires one to
\begin{enumerate}
    \item construct a complex functor $\cC_\pvV$ to serve as a ``candidate'' for $\PV$,
    \item prove that $\cC_\pvV$ is a lower bound for $\PV$, and
    \item prove that $\cC_\pvV$ is an upper bound for $\PV$.
\end{enumerate}
Notice that once (1) is accomplished, the main Galois connection (Theorem \ref{MAINTHEOREM}) means that establishing (2) is equivalent to establishing that all members of $\pvV$ are fixed by $\cC_\pvV$.
It should be noted that (3) is considerably more difficult to accomplish than (2)---and, moreover, that our framework does not seem to immediately make it any easier.
Hence our main concern from this point onward is to provide general tools for constructing ``candidate'' complex functors which satisfy (2).\footnotemark
\footnotetext{Note that since each lattice of $S$-complexes for a given finite semigroup $S$ is finite (and, moreover, each such lattice is computable), there are ``locally'' only finitely many possible candidates for $\cC(S)$ for any $\cC \in \CF$, all of which are---in and of themselves---computable.}

To this end, Section \ref{section:moduli} introduces \textit{moduli} (Definition \ref{mod:defn:moduli}), which are our basic tool for constructing complex functors. 
A modulus $\Lambda$ is a rule which assigns to each finite semigroup $S$ a (possibly empty) subset $\Lambda_S \subseteq \Po(S)$ in a manner which satisfies natural ``lift and push'' conditions with respect to morphisms.
The motivation for moduli comes from the observation that $\PV(S)$ is often characterized as being the minimal $S$-complex which is closed under unioning some distinguished subsets; e.g., subgroups for $\pvA$, $\GR$-classes for $\pvR$, $\GL$-classes for $\pvL$, and ``$\pvH$-kernels'' of subgroups for the pseudovariety of semigroups whose subgroups all belong to $\pvH$ (where $\pvH$ is some pseudovariety of groups).
Moduli generate complex functors---and hence complex monads as well---in a straightforward manner, and it will be shown that the fixed points of said complex functor are precisely those semigroups to whom said modulus assigns at most singletons (Theorem \ref{MODULARCONSTRUCTIONTHEOREM}).
Consequently, constructing a lower bound for $\PV$ is equivalent to defining a modulus which assigns at most singletons to members of $\pvV$.

Section \ref{section:effectiveness} discusses \textit{effective} moduli with respect to a pseudovariety $\pvV$, which are moduli whose induced complex monad is $\PV$. 
A number of known examples are provided to illustrate the utility of our language.

The next three sections are concerned with establishing transfer results for pointlikes along continuous operators on the lattice of pseudovarieties; which, as noted above, is the primary intended application of our framework.
Section \ref{section:transfer} introduces our approach, which is concerned with finding pairs of operators
    \[
        \lambda : \CF \arw \CF
            \quad \text{and} \quad
        \alpha : \PVAR \arw \PVAR
    \]
which make one or both of the diagrams
    \begin{equation*}
        \begin{tikzcd}
            \CF 
                \ar[rr, two heads, "\Fix"] 
                \ar[dd, "\lambda"'] 
        && 
            \PVAR
                \ar[dd, "\alpha"]
        && 
            \CF  
                \ar[dd, "\lambda"'] 
        && 
            \PVAR
                \ar[ll, hook', "\PLF"']
                \ar[dd, "\alpha"]
        \\&&& \text{and} &&&\\
            \CF
                \ar[rr, two heads, "\Fix"]
        && 
            \PVAR
        && 
            \CF
        &&
            \PVAR
                \ar[ll, hook', "\PLF"']
        \end{tikzcd}
    \end{equation*}
commute.
Given such a pair, we say that $\lambda$ satisfies the \textit{fixed point transfer condition} (respectively, \textit{pointlike transfer condition}) \textit{with respect to $\alpha$} if the left-hand (respectively, right-hand) diagram commutes (Definition \ref{tr:defn:conditions}).

To illustrate this language, Section \ref{section:reversal} provides an operator on $\CF$ which satisfies the pointlike transfer condition with respect to the operator $\pvV \mapsto \pvV^\rev$ induced by reversing the operations of members of $\pvV$.

Section \ref{section:contextspec} introduces \textit{context specifiers} (Definition \ref{cts:defn:contextspec}), which are moduli assigning sets of subsemigroups rather than of mere subsets.
Examples include the rules selecting ``$\pvW$-idempotent pointlike''---and, more generally, ``$\pvU$-like with respect to $\pvW$''---subsemigroups, subgroups, group kernels, local monoids, and the subsemigroups generated by idempotents or regular elements. 
Context specifiers induce join-preserving operators (and hence Galois connections) on $\PVAR$ which send a pseudovariety $\pvV$ to the pseudovariety consisting of semigroups for which all subsemigroups assigned by the context specifier belong to $\pvV$ (Definition \ref{cts:defn:contextgalois}). 
Many familar operators are shown to arise this way; examples include Mal'cev and generalized Mal'cev products---in one variable and of the forms $( - ) \malcev \pvW$ and $( - , \pvU ) \malcev \pvW$, respectively---as well as the ``subgroups belong to'', ``local monoids belong to'', ``idempotents generate'', and ``regular elements generate'' operators.

For each context specifier, we define an operator on $\CF$ which satisfies the fixed point transfer condition with respect to the aforementioned supremum-preserving operator on $\PV$ (Theorem \ref{cts:thm:moduli}).
The image of a decidable complex functor under one of these operators is decidable so long as it is induced by a decidable context specifier; this yields an ``automatic'' method of transferring lower bound results for pointlike sets along a large class of continuous operators on $\PVAR$.
This suggests a natural question: when do these operators satisfy the associated pointlike transfer condition as well?
This is likely to be a nontrivial question to answer in any generality, but it is a question which has considerable precedent in the literature.
A notable example is the main result of Steinberg and van Gool's paper \cite{PL-VARDETGRP-BEN}, which is equivalent to the statement that (the monad completion of) the operator on $\CF$ induced by the subgroup context satisfies the pointlike transfer condition with respect to the operator sending $\pvV \in \PVAR$ to the pseudovariety of semigroups whose subgroups all belong to $\pvV$.\footnotemark
\footnotetext{Their paper restricts its attention to the cases where $\pvV$ is a pseudovariety of groups, but these are easily seen to be the only ones which matter (consider $\pvV \cap \pvG$).}
Additionally, a conjecture of Steinberg in \cite{PL-JOINS-BEN} is equivalent to the statement that the operator on $\CF$ induced by the $\pvW$-idempotent pointlikes context satisfies the pointlike transfer condition with respect to $(-) \malcev \pvW$ whenever $\pvW$ is a pseudovariety of bands.

We conclude the paper in Section \ref{section:conclusion} with a discussion of possible future trajectories for the theory developed here. 

%

%% file: content/mainmatter/prelim.tex
%

\section{Preliminaries}      \label{section:prelim}

We assume that the reader is familiar with finite semigroup theory and has working knowledge of elementary category theory (including the standard limits and colimits, adjunctions, monads, and the Grothendieck construction) and lattice theory (particularly the theory of Galois connections).
For further background, see \cite{QTHEORY, ALM-BOOK, GRILLET} for semigroup theory; see \cite{CTCONTEXT, CATWORK, BARRCTCS} for category theory; and see \cite{COMPENDIUM} for lattice theory.

\subsection{Composition and application}
The majority of ``composable stuff'' in this paper will be composed from \textit{left to right}; this convention is sometimes called \textit{diagrammatic} order:
\begin{equation*}
            \begin{tikzcd}
            &
            \cdot 
            \arrow[dr, "b"]
            &
            \\
            \cdot 
            \arrow[ur, "a"] 
            \arrow[rr, "ab"]
            &&
            \cdot
            \end{tikzcd}
\end{equation*}
Morphisms between semigroups will always be composed this way.
Following this convention, the evaluation of a function $f$ at an element $x$ of its domain will be written as $(x)f$.

However, there are exceptions to this rule---generally involving functors (e.g. pointlikes)---in cases where there is a sufficiently powerful convention which demands otherwise.
In practice, these exceptions are made clear by context (along with generous use of brackets and parentheses). 
Moreover, both rule-abiding and rule-breaking situations will usually be accompanied by a diagram which makes the order in which things are to be composed obvious.

When it is useful to do so, the symbol $\fatsemi$ will be used to unambiguously indicate composition in diagrammatic order; e.g., $F \fatsemi G$ will always mean \textit{do $F$ then $G$}, regardless of any notational conventions regarding $F$ and $G$.

\subsection{Categorical notation}
Let $\sfC$ be a category.
We write $x \in \sfC$ to indicate that $x$ is an object of $\sfC$,
and the identity morphism at $x \in \sfC$ is denoted by $\id{x}$.
Moreover, given $x, y \in \sfC$, the set of $\sfC$-morphisms from $x$ to $y$ is denoted by $\sfC(x, y)$.

\subsection{Adjunctions and monads}
Recall that an \textbf{adjunction} is the data of a pair of functors $F : \sfC \rightarrow \sfD$ and $G : \sfD \rightarrow \sfC$ along with a natural isomorphism $\sfD( (-)F, -) \cong \sfC(-, (-)G)$.
This is represented in writing as $F \dashv G$ and diagrammatically as 
\begin{equation*}
            \begin{tikzcd}
            \sfC \arrow[rr, "F"', bend right=15, ""{name=B,above}] 
            & &
            \sfD. \arrow[ll, "G"', bend right=15, ""{name=A,below}].
            \ar[from=A, to=B, phantom, "\ddash"]
        \end{tikzcd}
\end{equation*}
It is at this point that the ``horizontal orientation'' conundrum must be confronted: we call $F$ the \textbf{lower adjoint} and $G$ the \textbf{upper adjoint}.

A \textbf{monad} on a category $\sfC$ is the data of a triple $(T, \varepsilon, \mu)$ consisting of
\begin{itemize}
    \item an endofunctor $T : \sfC \rightarrow \sfC$,
    \item a natural transformation $\varepsilon: \id{\sfC} \Rightarrow T$ called the \textbf{unit}, and
    \item a natural transformation $\mu : T^2 \Rightarrow T$ called the \textbf{multiplication}
\end{itemize}
such that both the \textit{unit diagram}
\begin{equation*}
            \begin{tikzcd}
                T 
                \ar[dr, Rightarrow, "\id{T}"'] 
                \ar[r, Rightarrow, "\varepsilon T"]
            &
                T^2   
                \ar[d, Rightarrow, "\mu"]
            &
                T
                \ar[l, Rightarrow, "T\varepsilon"']
                \ar[dl, Rightarrow, "\id{T}"]
            \\
            &
                T 
            &
            \end{tikzcd}
\end{equation*}
as well as the \textit{associativity diagram}
\begin{equation*}
            \begin{tikzcd}
                    T^3 
                    \arrow[d, Rightarrow, "\mu T"'] 
                    \arrow [rr, Rightarrow, "T \mu"]
                &&
                    T^2 
                    \arrow[d, Rightarrow, "\mu"]
                \\
                    T^2 
                    \arrow [rr, Rightarrow, "\mu"]
                && 
                    T
            \end{tikzcd}
\end{equation*}
commute.
In abstract terms, a monad is a monoid object in the category of endofunctors on $\sfC$.

\subsection{Order theoretic notation}
Let $(P , \leq)$ be a preorder.
A unique minimal upper bound for $X \subseteq P$ is called the \textbf{join} of $X$, and if such a thing exists it is denoted by $\bigvee X$; and, dually, a unique maximal lower bound for $X$ is called the \textbf{meet} of $X$ and is denoted by $\bigwedge X$ if such a thing exists.
If $X = \{ x, y \}$, then it is only polite to use infix notation for its join and meet, i.e., 
    \[
        x \vee y
        \quad \text{and} \quad
        x \wedge y
    \]
denote $\bigvee \{x , y\}$ and $\bigwedge \{x , y\}$, respectively.

\subsection{Complete lattices}
A poset wherein every subset has both a join and a meet is a \textbf{complete lattice}.
In a complete lattice, the meet and join are mutually determined by
    \begin{align*}
        \bigvee X &\;=\; \bigwedge \setst{\ell}{\textnormal{$x \leq \ell$ for all $x \in X$}}
        \\
        \bigwedge X &\;=\; \bigvee \setst{\ell}{\textnormal{$\ell \leq x$ for all $x \in X$}}
    \end{align*}
for every subset $X$ of said lattice.
The \textbf{top} of a complete lattice $L$ is the join of all of $L$, and the \textbf{bottom} of $L$ is the meet of all of $L$.

\subsection{Galois connections and closure operators}
A \textbf{Galois connection} between posets $P_1$ and $P_2$ is a pair of opposing monotone functions $\delta : P_1 \rightarrow P_2$ and $\gamma : P_2 \rightarrow P_1$ such that
    \[
        (x_1)\delta \leq x_2
        \qquad \Longleftrightarrow \qquad
        x_1 \leq (x_2)\gamma
    \]
for all $x_1 \in P_1$ and all $x_2 \in P_2$.

In other words, a Galois connection is an adjunction between posets, and the maps $\delta$ and $\gamma$ are the lower and upper adjoints, respectively.
Hence the notation $\delta \dashv \gamma$ will be used to indicate that a pair $(\delta , \gamma)$ of opposing maps between posets is a Galois connection.
Lower adjoints preserve all existing joins, and upper adjoints preserve all existing meets.

A \textbf{closure operator} on a poset $P$ is a map $c : P \rightarrow P$ which is monotone, increasing, and idempotent.

A closure operator $c$ on a complete lattice $L$ gives rise to a Galois connection in the following manner.
The set of fixed points of $c$---that is, elements $\ell \in L$ such that $(\ell)c = \ell$; or, equivalently, $\Img(c)$---is closed under meets.
Hence $\Img(c)$ is a complete lattice under the inherited meet and order whose join is determined, and there is a Galois connection
\begin{equation*}
            \begin{tikzcd}
            L \arrow[rr, two heads, "c"', shift right=2, ""{name=B,above}] 
            & &
            \Img(c) \arrow[ll, hook', "i"', shift right=2, ""{name=A,below}]
            \ar[from=A, to=B, phantom, "\ddash"]
        \end{tikzcd}
\end{equation*}
induced by the inclusion $\Img(c) \inj L$.

\subsection{Finite semigroups}
Let $S$ be a finite semigroup.
Let $S^I$ denote the semigroup obtained by adjoining an element $I$ and defining $x I = I x = x$ for every $x \in S^I$.
The set of idempotents of $S$ is denoted by $E(S)$ and the set of regular elements of $S$ is denoted by $\Reg(S)$.
Green's equivalence relations are denoted by $\GR$, $\GL$, $\GH$, and $\GJ$.

\subsection{The category of finite semigroups}
The category of finite semigroups and their morphisms---denoted by $\FSGP$---is \textit{Barr-exact}, i.e., it is a regular category in which congruences and kernel pairs coincide.\footnotemark
\footnotetext{See \cite{BARREXACT}.}
The empty semigroup $\varnothing$ and the trivial semigroup $\trivsgp$ are initial and terminal in $\FSGP$, respectively.

Monomorphisms in $\FSGP$ are precisely the injective morphisms and are denoted by arrows like $\inj$ with a hook.
Although not every epimorphism in $\FSGP$ is surjective,\footnotemark 
\footnotetext{See \cite[Example~3.1]{EPIDOMINION}.}
the \textit{regular} epimorphisms are precisely the surjective homomorphisms; these are denoted by arrows like $\surj$ with two heads.\footnotemark
\footnotetext{These arrow-writing conventions will carry over to other concrete categories as well.}

\subsection{Relational morphisms and division}
A \textbf{relational morphism} $\rho : S \rlm T$ is an equivalence class of spans in the category $\FSGP$ of the form
    \begin{equation*}
        \begin{tikzcd}
            \cdot 
            \arrow[r] 
            \arrow[d, two heads] 
        & 
            T 
        \\
            S 
        &
        \end{tikzcd}
    \end{equation*}
where the map to $S$ is a regular epimorphism, and where two such spans are equivalent if the natural maps from each apex to $S \times T$ have the same image.
This shared image is called the \textit{graph} of $\rho$ and is denoted by $\Gamma(S , \rho, T)$.
Individual spans in the equivalence class of $\rho : S  \rlm T$ are called \textbf{factorizations} of $\rho$.
The span
    \begin{equation*}
        \begin{tikzcd}
            \Gamma(S, \rho, T) 
            \arrow[r, "\rho_T"] 
            \arrow[d, two heads, "\rho_S"'] 
        & 
            T 
        \\
            S 
        &
        \end{tikzcd}
    \end{equation*}
where $\rho_S, \rho_T$ are the restricted natural projections is called the \textbf{canonical factorization} of $\rho : S \rlm T$.

By the \textbf{image} of $\rho : S \rlm T$, we will mean the image of the codomain projection $\rho_T : \Gamma(S, \rho, T) \rightarrow T$.
The image of $\rho$ will be denoted by $\Img(\rho)$.

Moreover, given $s_0 \in S$, let
    \[
        (s_0)\rho = \setst{t \in T}{(s_0 , t) \in \Gamma(S , \rho , T)};
    \]
and, similarly, let
    \[
        (t_0)\rho^\inv = \setst{s \in S}{(s , t_0) \in \Gamma(S, \rho, T)}
    \]
for any $t_0 \in T$.

A semigroup $S$ is said to \textbf{divide} a semigroup $T$ if $S$ is the homomorphic image of a subsemigroup of $T$---i.e., a subquotient. 
Equivalently, $S$ divides $T$ if and only if they are connected by some span of the form
    \begin{equation*}
        \begin{tikzcd}
            \cdot 
            \arrow[r, hook] 
            \arrow[d, two heads] 
        & 
            T 
        \\
            S 
        &
        \end{tikzcd}
    \end{equation*}
which is easily seen to be equivalent to the statement that there exists a relational morphism $\delta : S \rlm T$ for which $\delta_T : \Gamma(S, \delta, T) \inj T$ is a monomorphism.
In this case, we say that $\delta$ \textit{is} a \textbf{division}, and we write $\delta : S \sdiv T$.

Relational morphisms may be composed via pullback
    \begin{equation*}
        \begin{tikzcd}
            \cdot
            \ar[d, two heads]
            \ar[r]
            \ar[dd, two heads, bend right=35]
            \ar[rr, bend left=25]
            \ar[dr, phantom, "\lrcorner", very near start]
        &
            \cdot
            \ar[d, two heads]
            \ar[r]
        &
            U
        \\
            \cdot 
            \arrow[r] 
            \arrow[d, two heads] 
        & 
            T 
        &
        \\
            S 
        &
        \end{tikzcd}
    \end{equation*}
and the graph of the composition of $\rho : S \rlm T$ with $\mu : T \rlm U$ is given by
    \[
        \Gamma(S, \rho\mu , U) = \setst{(s , u)}{\textnormal{$(s , t_0) \in \Gamma(S , \rho , T)$ and $(t_0, u) \in \Gamma(T, \mu , U)$ for some $t_0 \in T$ }}.
    \]
Of course, morphisms are relational morphisms as well.
If $\varphi : S \rightarrow T$ is a morphism and $\mu : T \rlm U$ is a relational morphism, then
    \[
        \Gamma(S, \varphi\mu , U)
        = \setst{(s, u)}{(s\varphi , u) \in \Gamma(T, \mu, U)};
    \]
and if $\rho : S \rlm T$ is a relational morphism and $\psi : T \rightarrow U$ is a morphism, then
    \[
        \Gamma(S , \rho\psi, U) = \setst{(s, t\psi)}{(s , t) \in \Gamma(S , \rho , T)}.
    \]

\subsection{Pseudovarieties and continuous operators}
Recall that a \textbf{pseudovariety} is a class of finite semigroups which is closed under taking subobjects, homomorphic images, and finite products of its members.
In practice we identify a pseudovariety with the full subcategory of $\FSGP$ spanned by its members.

The collection of pseudovarieties is denoted by $\PVAR$; and if $\mathscr{A} \subseteq \FSGP$ is a class of finite semigroups, then $\PVGEN{\mathscr{A}}$ denotes the pseudovariety generated by semigroups in $\mathscr{A}$.

Crucially, $\PVAR$ is a complete lattice wherein the partial order is inclusion, the top and bottom are $\FSGP$ and $\pvtriv$ respectively, the meet is intersection, and the join is determined.
Alternatively, the join of $\pvV, \pvW \in \PVAR$ is given by
    \[
        \pvV \vee \pvW = \PVGEN{  \setst{V \times W}{V \in \pvV, \, W \in \pvW}  }.
    \]
Note that the join of a directed set of pseudovarieties is their union.

The lattice (see \cite{QTHEORY}) of continuous operators on the lattice of pseudovarieties is denoted by $\CntPV$.

\subsection{Common operators}
Given $\pvU, \pvV, \pvW \in \PVAR$, write
\begin{itemize}
    \item $\mathbb{E}\pvV = \setst{S}{\langle E(S) \rangle \in \pvV}$;
    \item $\mathbb{R}\pvV = \setst{S}{\langle\Reg(S)\rangle \in \pvV}$;
    \item $\mathbb{G}\pvV = \setst{S}{\text{all subgroups of $S$ belong to $\pvV$}}$;
    \item $\mathbb{L}\pvV = \setst{S}{\text{$eSe \in \pvV$ for all $e \in E(S)$}}$;
    \item $\pvV \ast \pvW = \PVGEN{ \setst{V \rtimes W}{\text{$V \in \pvV$ and $W \in \pvW$}}}$;
    \item $(\pvU, \pvV) \malcev \pvW$ for the pseudovariety of semigroups $S$ for which there exists a relational morphism $\rho : S \rlm W$ with $W \in \pvW$ such that, if $V$ is a subsemigroup of $W$ which belongs to $\pvV$, then $(V)\rho^\inv$ belongs to $\pvU$; and
    \item $\pvV \malcev \pvW = (\pvV, \pvtriv) \malcev \pvW $.
\end{itemize}

The last two operations are called the \textit{generalized Mal'cev product} and the \textit{Mal'cev product}, respectively.

%

%% file: content/mainmatter/complexes.tex
%

\section{Semigroup complexes}      \label{section:complexes}

\subsection{}
Recall the \textit{power functor} $\Po : \FSGP \rightarrow \FSGP$ which sends a finite semigroup $S$ to the semigroup $\Po(S)$ whose elements are the non-empty subsets of $S$ and which sends a morphism $\varphi : S \rightarrow T$ to the morphism
    \[
        \ext{\varphi} : \Po(S) \arw \Po(T)
            \quad \text{given by} \quad
        (X)\ext{\varphi} = \setst{(x)\varphi}{x \in X}.
    \]
The functor $\Po$ creates monomorphisms, regular epimorphisms, and isomorphisms.
Moreover, the collection of singleton embeddings
    \[
        \{-\} : S \longinj \Po(S)
            \quad \text{given by} \quad
        s \longmapsto \{s\}
    \]
as $S$ ranges over all finite semigroups are the components of a natural transformation from the identity functor on $\FSGP$ to $\Po$.
The semigroup of singletons of a finite semigroup $S$ is denoted by $\sing(S)$.

\begin{defn}    \label{com:defn:complex}
A (finite) \textbf{semigroup complex} is a pair $( S , \cK )$ consisting of a finite semigroup $S$ and a subsemigroup $\cK$ of $\Po(S)$ which
\begin{enumerate}
    \item contains $\sing(S)$ as a subsemigroup; and which 
    \item is closed under taking non-empty subsets of its members, meaning that if $X \in \cK$ then any $X_0 \in \Po(S)$ for which $X_0 \subseteq X$ is also a member of $\cK$.
\end{enumerate}
\end{defn}

\subsection{}
The name ``semigroup complexes'' comes from the fact that these are precisely the semigroup objects in the category of finite abstract simplicial complexes.
Hence given a semigroup complex $(S , \cK)$, we refer to $S$ as the \textit{vertex semigroup} and to $\cK$ as the \textit{face semigroup}.

\subsection{Morphisms}  \label{com:env:morphisms}
A morphism of semigroup complexes $\varphi : (S , \cK_S) \rightarrow (T, \cK_T)$ is given by a morphism $\varphi : S \rightarrow T$ in the category of finite semigroups with the property that $(X)\ext{\varphi} \in \cK_T$ for every $X \in \cK_S$.

Notice that $\varphi : S \rightarrow T$ satisfies this condition if and only if the restricted extension $\ext{\varphi} : \cK_S \rightarrow \cK_T$ is well-defined.

\begin{notation}
The category of semigroup complexes is denoted by $\COMCAT$.
\end{notation}

\subsection{}
The category $\COMCAT$ appears alongside a pair of functors and a natural transformation, which are arranged as in the diagram
\begin{equation*}
            \begin{tikzcd}
                \COMCAT
                    \arrow[rrr, bend left=15, ""{name=A,below}, "\comvert"]
                    \arrow[rrr, bend right=15, ""{name=B,above}, "\comface"']
            &
            &   
            &
                \FSGP   
            \ar[Rightarrow, from=A, to=B, "\:\comsingnat"]
        \end{tikzcd}
\end{equation*}
and which are defined as follows.
\begin{enumerate}
    \item The functor $\comvert$ is the evident forgetful functor which sends a semigroup complex to its vertex semigroup and which acts in the obvious unobtrusive manner on morphisms. 
    \item The functor $\comface$ sends semigroup complexes to their face semigroups and morphisms to the extensions of their vertex morphisms.
    \item The components of $\comsingnat$ are the singleton embeddings, i.e.,
        \[
            \comsingnat 
                \;\;=\;\;
            \Big( \{ - \} : S \longinj \cK \Big)_{(S, \cK) \,\in\, \COMCAT}.
        \]
\end{enumerate}
Alternatively---and more pleasantly---the diagram
\begin{equation*}
\begin{tikzcd}
&
            (S, \cK_S)
                \ar[ddddl, no head, dashed, bend right=25]
                    \ar[dddr, no head, dashed, bend left=25]
                \ar[rrrrrr, "\varphi"]
& & & 
& & & 
            (T, \cK_T)
                \ar[ddddl, no head, dashed, bend right=25]
                    \ar[dddr, no head, dashed, bend left=25]
&
\\
&&&&&&&&
\\
&&&&&&&&
\\
& & 
            \cK_S
                \ar[rrrrrr, "\ext{\varphi}" near start]
& & & 
& & & 
            \cK_T
\\
            S
                \ar[urr, hook, "\{-\}"]
                    \ar[rrrrrr, "\varphi"]
& & & 
& & &
            T 
                \ar[urr, hook, "\{-\}"]
& & &
\end{tikzcd}
\end{equation*}
illustrates the situation at a morphism in $\COMCAT$, where the left and right dashed lines represent $\comvert$ and $\comface$, respectively.

\subsection{Products and pullbacks}
The category $\COMCAT$ has all finite products; the product of semigroup complexes $(S_1, \cK_1)$ and $(S_2, \cK_2)$ is the semigroup complex
    \[
        (S_1, \cK_1) \times (S_2, \cK_2)
            \;\;=\;\;
        (S_1 \times S_2 ,\; \cK_1 \otimes \cK_2)
    \]
whose vertex semigroup is given by the product in $\FSGP$, and whose face semigroup is defined by
    \[
        \cK_1 \otimes \cK_2 
            \;\;=\;\;
        \setst{ Z \in \Po(S_1 \times S_2) }{ \text{ $(Z)\ext{\pi_{S_i}} \in \cK_i$ for $i = 1 , 2$} }.
    \]
Moreover, the map
    \[
        (- \times -) : \cK_1 \times \cK_2 \arw \cK_1 \otimes \cK_2
            \quad \text{given by} \quad
        (X_1, X_2) \longmapsto X_1 \times X_2
    \]
is a section of the natural map
    \[
        ( \ext{\pi_{S_1}} , \ext{\pi_{S_2}} ) : 
            \cK_1 \otimes \cK_2 \longsurj \cK_1 \times \cK_2
        \quad \text{given by} \quad
            Z \longmapsto \big( (Z)\ext{\pi_{S_1}} , (Z)\ext{\pi_{S_2}} \big).
    \]
The collections of maps
    \[
        \Big( (- \times -) : \cK_1 \times \cK_2 \longinj \cK_1 \otimes \cK_2 \Big)_{(S_1, \cK_1), (S_2 , \cK_2) \in \COMCAT}
    \]
and
    \[
        \Big( 
                ( \ext{\pi_{S_1}} , \ext{\pi_{S_2}} ) 
                : 
                \cK_1 \otimes \cK_2 \longsurj \cK_1 \times \cK_2 
        \Big)_{(S_1, \cK_1), (S_2 , \cK_2) \in \COMCAT}
    \]
are natural transformations indicated by the bold arrows in the diagram
\begin{equation*}
            \begin{tikzcd}
                    \COMCAT \times \COMCAT
                    \ar[dddd, "\comface \times \comface"'] 
                    \ar[rrr, "\times"]
                    \ar[ddddrrr, shift left=1, bend left=16, dashed, ""{name=A,below, near start}, ""{name=C,below, near end}]
                    \ar[ddddrrr, shift left=1, bend right=16, dashed, ""{name=B,above, near start}, ""{name=D,above, near end}]
                    \ar[Rightarrow, from=B, to=A]
                    \ar[Rightarrow, from=C, to=D]
                &&&
                    \COMCAT 
                    \ar[dddd, "\comface"]
                \\&&&\\&&&\\&&&\\
                    \FSGP \times \FSGP 
                    \ar[rrr, "\times"']
                &&& 
                    \FSGP.
            \end{tikzcd}
\end{equation*}

In fact, more is true: $\COMCAT$ has all finite pullbacks, which are computed as follows.
Given a cospan
\begin{equation*}
            \begin{tikzcd}
                &&
                    (S_1, \cK_1) 
                    \ar[dd, "\varphi_1"]
                \\ && \\
                    (S_2, \cK_2) 
                    \ar[rr, "\varphi_2"']
                && 
                    (T, \cL)
            \end{tikzcd}
\end{equation*}
in $\COMCAT$, the apex of the ensuing pullback is the semigroup complex
    \[
        (S_1, \cK_1) \times_{(T, \cL)} (S_2, \cK_2)
            \;\;=\;\;
        ( S_1 \times_T S_2 , \; \cK_1 \otimes_\cL \cK_2 ),
    \]
where $\cK_1 \otimes_\cL \cK_2$ is defined by
    \[
        \cK_1 \otimes_\cL \cK_2
            \;\;=\;\;
        \big( \cK_1 \otimes \cK_2 \big) \,\cap\, \Po(S_1 \times_T S_2).
    \]

\subsection{A local view of complexes}
Let $S$ be a finite semigroup.
A subsemigroup $\cK$ of $\Po(S)$ for which $(S, \cK)$ is a semigroup complex is called an \textbf{$S$-complex}. 

\begin{notation}
Let $S$ be a finite semigroup.
    \begin{enumerate}
    \item The set of $S$-complexes is denoted by $\Com{S}$, i.e.,
        \[
            \Com{S}
                =
            \setst{ \cK }{(S , \cK) \in \COMCAT }.
        \]
    \item Denote the closure of $\cX \subseteq \Po(S)$ under taking non-empty subsets of its members by $\ssc{\cX}$; i.e., 
        \[
            \ssc{\cX}
                =
            \setst{ Y \in \Po(S) }{ \text{$Y \subseteq X$ for some $X \in \cX$} }.
        \]
    \end{enumerate} 
\end{notation}

\subsection{Local lattices of complexes}    \label{com:env:scomlattice}
Given a finite semigroup $S$, the \textit{$S$-complex generation map} sending a set $\cX \subseteq \Po(S)$ to the $S$-complex
    \[
        \COMGEN{\cX}{S}
            =
        \ssc{ \langle \cX \cup \sing(S) \rangle }
    \]
is easily seen to be a closure operator on the power set lattice of $\Po(S)$ whose fixed points are precisely the $S$-complexes.

This closure operator induces a complete lattice structure on $\Com{S}$ in which
    \begin{enumerate}
    \item the partial order is inclusion;
    \item the top and bottom elements are $\Po(S)$ and $\sing(S)$, respectively;
    \item the join of $\cK_1, \cK_2 \in \Com{S}$ is given by
        \[
            \cK_1 \vee \cK_2 
                = 
            \COMGEN{\cK_1 \cup \cK_2}{S};
        \]
    \item and the meet of $\cK_1, \cK_2 \in \Com{S}$ is simply their intersection $\cK_1 \cap \cK_2$.
    \end{enumerate}
Note that the meet of $\cK_1, \cK_2 \in \Com{S}$ may alternatively be expressed as
    \[
        \cK_1 \cap \cK_2 
            = 
        \setst{ X_1 \cap X_2 }{ \text{ $X_i \in \cK_i$ for $i = 1, 2$ and $X_1 \cap X_2 \neq \varnothing$ } }
    \]
due to the downward closure condition.

\subsection{Underlying functors to lattices}    \label{com:env:underlyingfunctors}
A morphism $\varphi : S \rightarrow  T$ induces a Galois connection between the lattices of $S$- and $T$-complexes
\begin{equation*}
        \begin{tikzcd}
                \Com{S} 
                \arrow[rrr, "\psh{\varphi}"', bend right=18, ""{name=B,above}] 
            & & &
                \Com{T} 
                \arrow[lll, "\plb{\varphi}"', bend right=18, ""{name=A,below}]
            \ar[from=A, to=B, phantom, "\ddash"]
        \end{tikzcd}
\end{equation*}
wherein the lower adjoint is
    \[
            \psh{\varphi} : \Com{S} \arw \Com{T}
            \quad\quad \text{given by} \quad\quad
            (\cK_S)\psh{\varphi} = \COMGEN{\,\setst{(X)\ext{\varphi}}{X \in \cK_S}\,}{T}
    \]
and the upper adjoint is
    \[
        \plb{\varphi} : \Com{T} \arw \Com{S}
        \quad\quad \text{given by} \quad\quad
        (\cK_T)\plb{\varphi} = \setst{X \in \Po(S)}{(X)\ext{\varphi} \in \cK_T}.
    \]
Equipping the object map $\Com{(-)}$ with the action on morphisms sending $\varphi$ to the lower adjoint $\psh{\varphi}$ yields a covariant functor from $\FSGP$ to the category of complete finite lattices and join-preserving maps; 
moreover, equipping $\Com{(-)}$ with the action sending morphisms to their induced upper adjoints yields a contravariant functor from $\FSGP$ to the category of complete finite lattices and meet-preserving maps (which is equivalent to the opposite of its counterpart with join-preserving maps).

The category of semigroup complexes (alongside the functor $\comvert$) may be recovered from either of these functors via the Grothendieck construction.
Considering this fact alongside the observation in \ref{com:env:morphisms} yields the following proposition.

\begin{prop} \label{com:prop:morphismequiv}
Let $\varphi : S \rightarrow T$ be a morphism, and let $\cK_S$ and $\cK_T$ be $S$- and $T$-complexes, respectively.
The following are equivalent:
\begin{enumerate}
    \item The morphism $\varphi : (S, \cK_S) \rightarrow (T, \cK_T)$ exists in $\COMCAT$.
    \item $(\cK_S)\psh{\varphi} \subseteq \cK_T$. 
    \item $\cK_S \subseteq (\cK_T)\plb{\varphi}$.
    \item The extension $\ext{\varphi} : \cK_S \rightarrow \cK_T$ in $\FSGP$ is well-defined.
\end{enumerate}
\end{prop}

\subsection{}
Going forward, we will often be concerned with whether or not extensions of morphisms to certain complexes inherit standard properties from their originators.
The following results will be useful in these situations.

\begin{lem}     \label{com:lem:downsurj}
Let $\varphi : S \surj T$ be a regular epimorphism.
If $\cU_S$ and $\cU_T$ are subsemigroups of $\Po(S)$ and $\Po(T)$, respectively, for which the extension $\ext{\varphi} : \cU_S \surj \cU_T$ is (well-defined and) a regular epimorphism, then the same is true of $\ext{\varphi}  : \ssc{\cU_S} \surj \ssc{\cU_T}$.
\end{lem}

\begin{proof}
If $Y_0 \in \ssc{\cU_T}$ then there is some $Y \in \cU_T$ of which $Y$ is a subset. 
Ex hypothesi there is some $X \in \cU_S$ for which $(X)\ext{\varphi} = Y$, from which it follows that $(X \cap (Y_0)\varphi^\inv ) \ext{\varphi} = Y_0$.
Since $X \cap (Y_0)\varphi^\inv$ is a subset of $X$, the claim follows.
\end{proof}

\begin{prop} \label{com:prop:monoepiequiv}
Let $\varphi : (S, \cK_S) \rightarrow (T, \cK_T)$ be a morphism in $\COMCAT$.
\begin{enumerate}
    \item $\ext{\varphi} : \cK_S \rightarrow \cK_T$ is a monomorphism if and only if $\varphi : S \rightarrow T$ is a monomorphism.
    \item $\ext{\varphi} : \cK_S \rightarrow \cK_T$ is a regular epimorphism if and only if $\varphi : S \rightarrow T$ is a regular epimorphism and $(\cK_S)\psh{\varphi} = \cK_T$.
\end{enumerate}
\end{prop}

\begin{proof}
\hfill
\begin{enumerate}
\item The "if" direction follows from the fact that the restriction of an injection to any subset of its domain remains an injection.
Similarly, the "only if" direction follows from the observation that, if $\ext{\varphi} : \cK_S \inj \cK_T$ is a monomorphism, then commutivity of the diagram
\begin{equation*}
            \begin{tikzcd}
                    S 
                    \ar[dd, "\{-\}"', hook'] 
                    \ar[rr, "\varphi"]
                &&
                    T 
                    \ar[dd, "\{-\}", hook]
                \\&&\\
                    \cK_S 
                    \ar[rr, "\ext{\varphi}", hook]
                && 
                    \cK_T
            \end{tikzcd}
\end{equation*}
forces $\varphi$ to be a monomorphism.
\item Begin with the ``if'' direction.
Notice that if $\varphi$ is a regular epimorphism then the image of $\ext{\varphi} : \cK_S \rightarrow \Po(T)$ contains the singletons of $T$.
Hence
    \begin{align*}
        (\cK_S)\psh{\varphi} 
            &=
        \COMGEN{ \Img(\ext{\varphi}) }{T}   \\
            &=
        \ssc{ \langle \Img(\ext{\varphi}) \cup \sing(T) \rangle }   \\
            &=
        \ssc{ \Img(\ext{\varphi}) }.
    \end{align*}
It follows from Lemma \ref{com:lem:downsurj} that $\ext{\varphi} : \ssc{\cK_S} \surj (\cK_S)\psh{\varphi}$ is a regular epimorphism; which, since $\ssc{\cK_S} =\cK_S$, yields the desired property. 

The ``only if'' direction is a consequence of the minimality of $(\cK_S)\psh{\varphi}$ amongst $T$-complexes containing the image of $\ext{\varphi} : \cK_S \rightarrow \Po(T)$.
\end{enumerate}
\end{proof}

%

%% file: content/mainmatter/complexfunctors.tex
%

\section{Complex functors}      \label{section:complexfunctors}

\begin{defn}    \label{cf:defn:complexfunctor}
A \textbf{complex functor} is an endofunctor $\cC : \FSGP \rightarrow \FSGP$ which
\begin{enumerate}
    \item sends a finite semigroup $S$ to an $S$-complex $\cC(S)$,
    \item acts on morphisms by sending $\varphi : S \rightarrow T$ to its (restricted) extension $\ext{\varphi} : \cC(S) \rightarrow \cC(T)$, and which
    \item preserves regular epimorphisms.
\end{enumerate}
\end{defn}

\begin{notation}
Let $\CF$ denote the collection of complex functors.
\end{notation}

\subsection{}       \label{cf:env:cflift}
Equivalently, a complex functor $\cC$ is an endofunctor on $\FSGP$ for which there exists a functor
    \[
        \comlift{\cC} : \FSGP \arw \COMCAT
    \]
which is a lift of $\cC$ through the face functor $\comface$---in the sense that the diagram
    \begin{equation*}
        \begin{tikzcd}
            &&
                \COMCAT 
                    \arrow[dd, "\comface"]
            \\&&\\
                \FSGP 
                    \arrow[uurr, "\comlift{\cC}"] 
                    \arrow[rr, "\cC"']
            &&
                \FSGP
        \end{tikzcd}
    \end{equation*}
commutes---and which satisfies the following conditions.
\begin{enumerate}
    \item The diagram of functors
    \begin{equation*}
        \begin{tikzcd}
            &&
                \COMCAT 
                    \arrow[dd, "\comvert"]
            \\&&\\
                \FSGP 
                    \arrow[uurr, "\comlift{\cC}"] 
                    \arrow[rr, "\id{\FSGP}"']
            &&
                \FSGP
        \end{tikzcd}
    \end{equation*}
    commutes; i.e., the action of $\comlift{\cC}$ is of the form
    \begin{equation*}
        \begin{tikzcd}
            S
            \ar[dd, "\varphi"', ""{name=A,above}]
        & & & 
                (S, \cC(S))
                \ar[dd, "\varphi", ""{name=B,below}]
            \ar[from=A, to=B, phantom, "{\xmapsto{ \quad \comlift{\cC} \quad }}"]
        \\ &&& \\
            T        
        & & &           
                (T, \cC(T))
        \end{tikzcd}
    \end{equation*}
    for any morphism $\varphi : S \rightarrow T$.
    \item The functor $\comlift{\cC}$ sends a regular epimorphism $\varphi : S \surj T$ to the morphism
        \[
            \varphi : (S, \cC(S)) \arw (T, (\cC(S))\psh{\varphi})
        \]
    in the category $\COMCAT$.
\end{enumerate}
Note that the second condition on $\comlift{\cC}$ is equivalent to the condition that $\cC$ preserves regular epimorphisms by Proposition \ref{com:prop:monoepiequiv}.

\begin{prop}    \label{cf:prop:createmor}
Complex functors create monomorphisms, regular epimorphisms, and isomorphisms.
\end{prop}

\begin{proof}
This follows immediately from Proposition \ref{com:prop:monoepiequiv}.
\end{proof}

\subsection{The singleton transformation}   \label{cf:env:singnat}
If $\cC$ is a complex functor, then there is a natural transformation $\sigma^\cC : \id{\FSGP} \Rightarrow \cC$ whose components are the singleton embeddings; that is,
    \[
        \sigma^\cC_S = \{-\} : S \longinj \cC(S)
    \]
for every finite semigroup $S$.

Moreover---building on \ref{cf:env:cflift}---the natural transformation $\sigma^\cC$ may be obtained via whiskering as illustrated in the diagram
\begin{equation*}
\begin{tikzcd}
    \COMCAT
        \ar[rrrr, bend left=15, "\comvert", ""{name=U, below}]
        \ar[rrrr, bend right=15, "\comface"', ""{name=D} ]
        \ar[Rightarrow, from=U, to=D, "\comsingnat"]
&
&&
&
    \FSGP
        \ar[ddd, equal]
\\ &&&& \\ &&&& \\
    \FSGP
        \ar[uuu, "\comlift{\cC}"]
            \ar[rrrr, bend left=15, "\id{\FSGP}", ""{name=A, below}]
            \ar[rrrr, bend right=15, "\cC"', ""{name=B}]
        \ar[Rightarrow, from=A, to=B, "\sigma^\cC"]
&
&&
&
    \FSGP.
\end{tikzcd}
\end{equation*}

\subsection{The lattice of complex functors}        \label{cf:env:cflattice}
There is a natural partial order on $\CF$ which is inherited pointwise from the various lattices $\Com{S}$ as $S$ ranges over all finite semigroups; concretely, this order is given by
    \[
        \cC_1 \leq \cC_2 
            \quad \Longleftrightarrow \quad
        \text{$\cC_1(S) \subseteq \cC_2(S)$ for all $S \in \FSGP$}
    \]
for all $\cC_1, \cC_2 \in \CF$.
The functors $\Po$ and $\sing$ are respectively maximal and minimal in $\CF$ with respect to this order.

Moreover, it is straightforward to show that the pointwise join
    \[
        \cC_1 \vee \cC_1
            \quad \text{given by} \quad
        [\cC_1 \vee \cC_2](S) = \cC_1(S) \vee \cC_2(S)
    \]
is the join of $\cC_1$ and $\cC_2$ in $\CF$.
This establishes that $\CF$ has all finite joins, and a routine argument via ``locally'' reducing infinite joins to finite joins at a given finite semigroup yields the existence of arbitrary joins in $\CF$ as well.

Consequently, $\CF$ is a complete lattice under the pointwise order and join, with the pointwise top and bottom, and whose meet is determined in general.

\begin{rmk}     \label{cf:rmk:ordernattr}
The partial order on $\CF$ may be realized ``concretely'' as natural transformations in two ways:
given $\cC_1 , \cC_2 \in \CF$ with $\cC_1 \leq \cC_2$, the components
    \[
        \Big( \id{S} : (S, \cC_1(S)) \arw (S, \cC_2(S)) \Big)_{S \in \FSGP}
    \]
constitute a natural transformation $\comlift{\cC_1} \Rightarrow \comlift{\cC_2}$ which, as shown in the diagram
\begin{equation*}
\begin{tikzcd}
    \FSGP
        \ar[ddd, equal]
        \ar[rrrr, bend left=15, "\comlift{\cC_1}", ""{name=U, below}]
        \ar[rrrr, bend right=15, "\comlift{\cC_2}"', ""{name=D} ]
        \ar[Rightarrow, from=U, to=D]
&
&&
&
    \COMCAT
        \ar[ddd, "\comface"]
\\ &&&& \\ &&&& \\
    \FSGP
        \ar[rrrr, bend left=15, "\cC_1", ""{name=A, below}]
        \ar[rrrr, bend right=15, "\cC_2"', ""{name=B}]
        \ar[Rightarrow, from=A, to=B]
&
&&
&
    \FSGP.
\end{tikzcd}
\end{equation*}
induces via whiskering a natural transformation $\cC_1 \Rightarrow \cC_2$ whose component at a given $S \in \FSGP$ is the inclusion $\cC_1(S) \inj \cC_2(S)$.
\end{rmk}

\begin{rmk}
The failure of meets in $\CF$ to be pointwise in general is due to the possible failure of pointwise meets to preserve regular epimorphisms.
In concrete terms, if $\cC_1$ and $\cC_2$ are complex functors and $\varphi : S \surj T$ is a regular epimorphism, then the extension
    \[
        \ext{\varphi} : \cC_1(S) \cap \cC_2(S) \arw \cC_1(T) \cap \cC_2(T)
    \]
might fail to be a regular epimorphism (although it is always well-defined).
This is related to the fact that the map $\psh{\varphi} : \Com{S} \rightarrow \Com{T}$ does not necessarily preserve meets (as it is a lower adjoint).

However, there are cases where the meet is pointwise; the following proposition provides an example of this.
\end{rmk}

\begin{prop} \label{cf:prop:filtdirect}
If $\{ \cC_\alpha \}_{\alpha \in I} \subseteq \CF$ is filtered, then
    \[
        \left[ \bigwedge_\alpha \cC_\alpha \right] (S) 
            =
        \bigcap_\alpha \cC_\alpha(S);
    \]
and, dually, if $\{ \cC_\alpha \}_{\alpha \in I} \subseteq \CF$ is directed, then 
    \[
        \left[ \bigvee_\alpha \cC_\alpha \right] (S) 
            =
        \bigcup_\alpha \cC_\alpha(S)
    \]
for all finite semigroups $S$.
\end{prop}

\begin{proof}
For any finite semigroup $S$, the set $\{ \cC_\alpha(S) \}_{\alpha \in I}$ is finite, and is directed or filtered if $\{ \cC_\alpha \}_{\alpha \in I}$ is directed or filtered in $\CF$.
The proposition then follows from the fact that any directed or filtered subset of a finite lattice has an upper or lower bound, respectively, which belongs to the set.
\end{proof}

\begin{defn}
A finite semigroup $S$ is said to be a \textbf{fixed point} of a complex functor $\cC$ if $\cC(S) = \sing(S)$.
\end{defn}

\subsection{}
Equivalently, $S \in \FSGP$ is a fixed point of $\cC \in \CF$ if and only if its component $\sigma^\cC_S : S \iso \cC(S)$ is an isomorphism (recall \ref{cf:env:singnat}). 

\begin{notation}
The set of fixed points of $\cC \in \CF$ is denoted by $\Fix(\cC)$; i.e.,
    \[
        \Fix(\cC)
            =
        \setst{ S \in \FSGP }{ \cC(S) = \sing(S) }.
    \]
\end{notation}

\begin{prop} \label{FIRSTPUNCHLINE}
\hfill
\begin{enumerate}
    \item If $\cC$ is a complex functor, then $\Fix(\cC)$ is a pseudovariety.
    \item The map $\Fix : \CF \rightarrow \PVAR$ is antitone and takes joins in $\CF$ to intersections in $\PVAR$.
    Moreover, $\Fix(\Po) = \pvtriv$ and $\Fix(\sing) = \FSGP$.
    \item Consequently, $\Fix$ has an antitone upper adjoint
        \[
            \MAXCF : \PVAR^\op \arw \CF 
                \quad \text{given by} \quad
            \MAXCF[\pvV] 
                =
            \bigvee \setst{ \cC \in \CF }{ \pvV \subseteq \Fix ( \cC ) },
        \]
    thus yielding an antitone Galois connection
    \begin{equation*}
            \begin{tikzcd}
            \CF \arrow[rrr, "\Fix"', bend right=15, ""{name=B,above}] 
            & & &
            \PVAR^\op \arrow[lll, "\MAXCF"', bend right=15, ""{name=A,below}]
            \ar[from=A, to=B, phantom, "\ddash"]
        \end{tikzcd}
    \end{equation*}
    between the complete lattices of complex functors and pseudovarieties.
\end{enumerate}
\end{prop}

\begin{proof}
\hfill
\begin{enumerate}
    \item We must show that $\Fix(\cC)$ is closed under subobjects, homomorphic images, and binary products.
    
    First, if $\varphi: U \inj S$ is a monomorphism and $S$ is a fixed point of $\cC$, then preservation of monomorphisms forces $U$ to be a fixed point of $\cC$ as well.
    
    Next, if $S$ is again a fixed point of $\cC$ and $\varphi : S \surj T$ is a regular epimorphism, then the second condition discussed in \ref{cf:env:cflift} implies that
        \[
            \cC(T) 
                = 
            (\cC(S))\psh{\varphi}
                =
            (\sing(S))\psh{\varphi},
        \]
    which in turn implies that $\cC(T) = \sing(T)$ since $\psh{\varphi}$ preserves joins. 
    
    Finally, if $S$ and $T$ are both fixed points of $\cC$, then the image of any $Z \in \cC(S \times T)$ under the extension of either projection is a singleton, and hence $Z$ must be a singleton as well.
    Hence $S \times T$ is fixed by $\cC$ as well, and the claim follows.
    \item Antitonicity of $\Fix$ is easily verified.
    
    Regarding joins, note that since the action of an arbitrary join in $\CF$ on a given finite semigroup may be reduced to the action of some finite join, it suffices to show that binary joins are taken to binary intersections---a fact which is easily verified.
    
    Finally, it is obvious that all finite semigroups are fixed by $\sing$; and the observation that the order of $\Po(S)$ is strictly greater than the order of $S$ if and only if $S$ is non-trivial yields the remaining claim. 
    \item This follows immediately from claims (1) and (2).
\end{enumerate}
\end{proof}

\begin{rmk}
This mysterious upper adjoint will turn out to be the map sending pseudovarieties to their respective pointlike functors.
\end{rmk}

%

%% file: content/mainmatter/complexmonads.tex
%

\section{Complex monads}      \label{section:complexmonads}

Recall that $\Po$ may be equipped with the structure of a monad $(\Po, \sigma, \mu)$, where the components of $\sigma : \id{\FSGP} \Rightarrow \Po$ are the singleton embeddings, and the components of $\mu : \Po^2 \Rightarrow \Po$ are the union maps, i.e., 
    \[
        \mu_S 
            =
        \bigcup (-) : \Po^2(S) \longsurj \Po(S)
    \]
for every finite semigroup $S$.

\begin{defn}    \label{cm:defn:complexmonad}
A \textbf{complex monad} is a complex functor $\cC$ which may be equipped with the additional structure required to yield a submonad of $(\Po, \sigma, \mu)$. 

Concretely, this means that $\cC$ admits a monad structure $(\cC , \sigma^\cC , \mu^\cC)$ where
\begin{enumerate}
    \item the components of $\sigma^\cC : \id{\FSGP} \Rightarrow \cC$ are the singleton embeddings, and
    \item the components of $\mu^\cC : \cC^2 \Rightarrow \cC$ are the union maps, i.e.,
    \[
        \mu^\cC_S 
            =
        \bigcup (-) : \cC^2(S) \longsurj \cC(S)
    \]
    for every finite semigroup $S$.
\end{enumerate}
In practice, we will omit specific reference to the natural transformations (since they are of a specified form) and refer to the complex functor $\cC$ and the monad of which it is a part interchangeably as $\cC$ whenever context allows.
\end{defn}

\begin{notation}
The collection of complex monads is denoted by $\CM$.
\end{notation}

\begin{lem} \label{cm:lem:unionnat}
If $\cC_1$, $\cC_2$, and $\cC_3$ are complex functors satisfying
    \[
        \cX \in \cC_1(\cC_2(S)) 
            \quad \Longrightarrow \quad 
        \bigcup \cX \in \cC_3(S)
    \]
for every finite semigroup $S$, then the induced morphisms
    \[
        \left(\; \bigcup (-) : \cC_1(\cC_2(S)) \arw \cC_3(S) \;\right)_{S \in \FSGP}
    \]
are natural.
Moreover, if $\cC_1$ and $\cC_2$ are complex functors which satisfy
    \[
        \cX \in \cC_1(\cC_2(S)) 
            \quad \Longrightarrow \quad 
        \bigcup \cX \in \cC_2(S)
    \]
for every finite semigroup $S$, then the components of the induced natural transformation are all regular epimorphisms.
\end{lem}

\begin{proof}
For the first claim, it is sufficient to show that
\begin{equation*}
            \begin{tikzcd}
                    \cC_1(\cC_2(S)) 
                    \arrow[dd, "\bigcup(-)"'] 
                    \arrow [rr, "\ext{\ext{\varphi}}"]
                &&
                    \cC_1(\cC_2(T))
                    \arrow[dd, "\bigcup(-)"]
                \\
                &&
                \\
                    \cC_3(S) 
                    \arrow [rr, "\ext{\varphi}"]
                && 
                    \cC_3(T)
            \end{tikzcd}
\end{equation*}
is commutative for every $\varphi : S \rightarrow T$.
Given $\cX \in \cC_1(\cC_2(S))$, the calculation
    \begin{align*}
        \bigcup \left[ (\cX)\ext{\ext{\varphi}} \right] 
        &\;\;=\;\;
        \bigcup \setst{(X)\ext{\varphi}}{X \in \cX}&
        \\
        &\;\;=\;\;
        \setst{(x)\varphi}{\textnormal{$x \in X$ for some $X \in \cX$}}&
        \\
        &\;\;=\;\; \left( \bigcup \cX \right) \ext{\varphi}
    \end{align*}
establishes the desired commutativity.
The second claim follows immediately from the fact that $\sing(\cC_2(S))$ is a subsemigroup of $\cC_1(\cC_2(S))$.
\end{proof}

\begin{prop}    \label{cm:prop:monadcondition}
If $\cC$ is a complex functor which satisfies
    \[
        \cX \in \cC^2(S) 
            \quad \Longrightarrow \quad 
        \bigcup \cX \in \cC(S)
    \]
for every finite semigroup $S$, then $\cC$ is a complex monad.
\end{prop}

\begin{proof}
Naturality of the singleton embeddings was established in \ref{cf:env:singnat}, and naturality of the union maps is a particular case of Lemma \ref{cm:lem:unionnat}.
Let $\sigma$ and $\mu$ denote the singleton and union transformations with respect to $\cC$, respectively.
There are two diagrams whose commutativity must be verified.

The first of these is the ``unit'' diagram:
\begin{equation*}
            \begin{tikzcd}
                \cC 
                \ar[dr, Rightarrow, "\id{\cC}"'] 
                \ar[r, Rightarrow, "\sigma \cC"]
            &
                \cC^2   
                \ar[d, Rightarrow, "\mu"]
            &
                \cC
                \ar[l, Rightarrow, "\cC\sigma"']
                \ar[dl, Rightarrow, "\id{\cC}"]
            \\
            &
                \cC 
            &
            \end{tikzcd}
\end{equation*}
The components of the natural transformation $\sigma \cC$ are the maps
    \[
        (\sigma\cC)_S = \sigma_{\cC(S)} : \cC(S) \longinj \cC^2(S)
        \quad\;\; \textnormal{given by} \quad\;\;
        X \longmapsto \{ X \},
    \]
and those of $\cC \sigma$ are 
    \[
        (\cC\sigma)_S = \ext{\sigma_{S}} : \cC(S) \longinj \cC^2(S)
        \quad\;\; \textnormal{given by} \quad\;\;
        X \longmapsto \setst{\{x\}}{x \in X}
    \]
for each finite semigroup $S$.
The observation that
    \[
        \bigcup \setst{\{x\}}{x \in X} 
        \;\; = \;\; 
        \bigcup \{X\} 
        \;\; = \;\;
        X
    \]
for any $X \in \cC(S)$ establishes the desired commutativity.

Next in line is the ``associativity'' diagram:
\begin{equation*}
            \begin{tikzcd}
                    \cC^3 
                    \arrow[d, Rightarrow, "\mu \cC"'] 
                    \arrow [rr, Rightarrow, "\cC \mu"]
                &&
                    \cC^2 
                    \arrow[d, Rightarrow, "\mu"]
                \\
                    \cC^2 
                    \arrow [rr, Rightarrow, "\mu"]
                && 
                    \cC
            \end{tikzcd}
\end{equation*}
The components of $\mu \cC$ are the maps
    \[
        (\mu \cC)_S = \mu_{\cC(S)} : \cC^3(S) \longsurj \cC^2(S)
        \;\; \textnormal{given by} \;\;
        \Sigma \longmapsto \setst{X}{\textnormal{$X \in \cX$ for some $\cX \in \Sigma$}}
    \]
and the components of $\cC \mu$ are 
    \[
        (\cC \mu)_S = \ext{\mu_{S}} : \cC^3(S) \longsurj \cC^2(S)
        \quad \textnormal{given by} \quad
        \Sigma \longmapsto \setst{\bigcup \cX}{\cX \in \Sigma}
    \]
for each finite semigroup $S$.

Now, let $\Sigma \in \cC^3(S)$.
Following $\Sigma$ through $(\mu \cC)_S$ yields
    \begin{align*}
        \bigcup \left[ (\Sigma)(\mu \cC)_S \right]
        &\;\;=\;\; \bigcup \setst{X}{\textnormal{$X \in \cX$ for some $\cX \in \Sigma$}} &
        \\
        &\;\;=\;\; \setst{x}{\textnormal{$x \in X$ for some $X \in \cX \in \Sigma$}};&
    \end{align*}
while following $\Sigma$ through $(\cC \mu)_S$ yields
    \begin{align*}
        \bigcup \left[ (\Sigma)(\cC \mu)_S \right]
        &\;=\;\; \bigcup \setst{\bigcup \cX}{\cX \in \Sigma}&
        \\
        &\;=\;\; \bigcup \setst{\setst{x}{\textnormal{$x \in X$ for some $X \in \cX$}}}{\cX \in \Sigma}&
        \\
        &\;=\;\; \setst{x}{\textnormal{$x \in X$ for some $X \in \cX \in \Sigma$}},&
    \end{align*}
from which the proposition follows.
\end{proof}

\subsection{Iterative union closure}
Given $\cC \in \CF$, define an increasing sequence
    \[
       \cflevel{\cC}{0}
            \,\leq\,
        \cflevel{\cC}{1}
            \,\leq\,
        \cflevel{\cC}{2}
            \,\leq\,
            \cdots
            \,\leq\,
        \cflevel{\cC}{n}
            \,\leq\,
        \cflevel{\cC}{n+1}
            \,\leq\, 
            \cdots
    \]
of complex functors recursively by setting
    \[
        \cflevel{\cC}{0}(S) = \sing(S)
            \quad \text{and} \quad
        \cflevel{\cC}{n+1}(S) = 
                    \ssc{ \setst{\bigcup \cA }{ \cA \in \cC \Big( \cflevel{\cC}{n}(S) \Big) } }
    \]
for every finite semigroup $S$.

Evaluating this recursion scheme at a chosen $S$ yields the diagram
\begin{equation*}
    \begin{tikzcd}
        \cC \big(\, \cflevel{\cC}{0}(S) \,\big)
            \ar[ddrr, "\bigcup(-)"]     
    &&
        \cC \big(\, \cflevel{\cC}{1}(S) \,\big) 
            \ar[ddrr, "\bigcup(-)"]   
    &&
        \cC \big(\, \cflevel{\cC}{2}(S) \,\big) 
            \ar[ddrr, phantom, "\ddots", shift left=3]          
    &&
        \cdots
    \\&&&&&&\\
        \cflevel{\cC}{0}(S)  
            \ar[rr, hook, "i"'] 
            \ar[uu, hook, "\{-\}"' near start]   
    &&
        \cflevel{\cC}{1}(S) 
            \ar[rr, hook, "i"'] 
            \ar[uu, hook, "\{-\}"' near start]  
    &&
        \cflevel{\cC}{2}(S)
            \ar[rr, hook, "i"'] 
            \ar[uu, hook, "\{-\}"' near start] 
    &&
        \cdots
    \end{tikzcd}
\end{equation*}
wherein every visible arrow is natural: the vertical arrows by \ref{cf:env:singnat}, the horizontal arrows by Remark \ref{cf:rmk:ordernattr}, and the diagonal arrows by Lemma \ref{cm:lem:unionnat}.

\begin{rmk}
It is readily verified that $\cflevel{\cC}{1} = \cC$ for any $\cC \in \CF$.
\end{rmk}

\begin{lem} \label{cm:lem:iterateterminate}
Let $\cC$ be a complex functor.
For every finite semigroup $S$ there is a natural number $N$ such that
    \[
        \cflevel{\cC}{N}(S) \;=\; \cflevel{\cC}{N+k}(S)
            \qquad
        \text{for all $k \in \mathbb{N}$}.
    \]
In this situation, $\cflevel{\cC}{N}(S)$ is a retract of $\cC \big( \cflevel{\cC}{N}(S) \big)$ as displayed in the diagram 
\begin{equation*}
            \begin{tikzcd}
            \cflevel{\cC}{N}(S) \arrow[rr, hook', "\{-\}"', shift right=2] 
            & &
            \cC \big( \cflevel{\cC}{N}(S) \big). \arrow[ll, two heads, "\bigcup(-)"', shift right=2]
        \end{tikzcd}
\end{equation*}
\end{lem}

\begin{proof}
Evaluation at $S$ yields an ascending chain 
    \[
        \cflevel{\cC}{0}(S) 
            \;\;\subseteq\;\; 
        \cflevel{\cC}{1}(S) 
            \;\;\subseteq\;\; \cdots \;\;\subseteq\;\; 
        \cflevel{\cC}{n}(S)
            \;\;\subseteq\;\; 
        \cflevel{\cC}{n+1}(S) 
            \;\;\subseteq\;\; \cdots
    \]
in the lattice $\Com{S}$, which---by finiteness---is guaranteed to converge.
The rest follows via an easy application of Lemma \ref{cm:lem:unionnat}.
\end{proof}

\begin{defn}    \label{cm:defn:moncom}
The \textbf{monad completion} of $\cC \in \CF$ is given by
    \[
        \moncom{\cC}(S) \;=\; \varinjlim \cflevel{\cC}{\bullet}(S)
    \]
for each finite semigroup $S$.
\end{defn}

\begin{prop} \label{SECONDPUNCHLINE}
\hfill
\begin{enumerate}
    \item If $\cC$ is a complex functor, then $\moncom{\cC}$ is a complex monad.
    \item The monad completion map
        \[
            \moncom{(-)} : \CF \arw \CF
        \]
    is a closure operator whose image is $\CM$.
    Consequently, $\CM$ is a complete lattice under the standard inherited and induced data.
    \item The fixed points of a complex functor $\cC$ and its monad completion $\moncom{\cC}$ coincide.
    \item Altogether, these claims yield a commuting diagram
    \begin{equation*}
    \begin{tikzcd}
            \CF
            \ar[rrrr, shift right=1, bend left=15, ""{name=B,above}, two heads, "\moncom{(-)}"']
            \ar[dddrr, shift left=1, bend right=15, ""{name=C}, "\Fix"]
            &       &       &       &
            \CM
            \ar[llll, shift right=3, bend right=15, ""{name=A,below}, hook', "i"']
            \ar[dddll, shift right=1.5, bend left=15, ""{name=E}, "\Fix"']
            \ar[from=A, to=B, phantom, "\scriptstyle{\ddash}"]
        \\
            &       &       &       &
        \\
            &       &       &       &
        \\
            &       &
            \PVAR^\op
            \ar[uuull, shift left=3, bend left=15, ""{name=D}, "\MAXCF"]
            \ar[uuurr, shift right=3, bend right=15, ""{name=F}, "\MAXCF"']
            &       &
        \ar[from=C, to=D, phantom, "\scriptstyle{\ddash}" rotate=135]
        \ar[from=E, to=F, phantom, "\scriptstyle{\ddash}" rotate=-135]
    \end{tikzcd}
    \end{equation*}
    of Galois connections between complete lattices.
\end{enumerate}
\end{prop}

\begin{proof}
\hfill
\begin{enumerate}
    \item Let $S$ be a finite semigroup.
    We will show that 
    \begin{equation} \label{cm:unionproperty}
        \Sigma \in \cflevel{\cC}{n} \left( \moncom{\cC}(S) \right)
        \quad \Longrightarrow \quad 
        \bigcup \Sigma \in \moncom{\cC}(S) \tag{$\star$}
    \end{equation}
    holds for all $n$.
    This will yield the claim by way of Proposition \ref{cm:prop:monadcondition}.
    Before we begin, note that if (\ref{cm:unionproperty}) holds for some $n$ then the union map
        \[
            \bigcup (-) \,:\, \cflevel{\cC}{n} \left( \moncom{\cC}(S) \right) \longsurj \moncom{\cC}(S)
        \]
    is a regular epimorphism.

    Proceed by induction on $n$.
    The case where $n=0$ is obvious, and the case where $n=1$ is an immediate consequence of Lemma \ref{cm:lem:iterateterminate}.

    Suppose that (\ref{cm:unionproperty}) holds for $n$, and let
        \[
            \mu \;=\; 
                \bigcup (-) \,:\, 
            \cflevel{\cC}{n} \left( \moncom{\cC}(S) \right) \longsurj \moncom{\cC}(S)
        \]
    denote the hypothesized map.
    As should be routine by now, the extension
        \[
            \ext{\mu} : \cC \left[ \cflevel{\cC}{n} \left( \moncom{\cC}(S) \right) \right] 
            \surj \cC \left[ \moncom{\cC}(S) \right]
                \;\;\; \textnormal{by} \;\;\;
            (\Sigma)\ext{\mu} = \setst{\bigcup \cA}{\cA \in \Sigma}
        \]
    is a regular epimorphism as well.
    Now, set
        \[
            \scK 
                \:=\:
            \setst{
                \bigcup \Sigma}{
                \Sigma \in 
                \cC \left[ \cflevel{\cC}{n} \left( \moncom{\cC}(S) \right) \right]
            },
    \]
    i.e., $\scK$ is the image of the map
        \[
            \bigcup (-) : 
                \cC \left[ \cflevel{\cC}{n} \left( \moncom{\cC}(S) \right) \right] 
                    \longsurj 
                \scK;
        \]
    and moreover $\scK$ has the property that $\ssc{\scK} = \cflevel{\cC}{n+1} \left( \moncom{\cC}(S) \right)$.

    The same argument used to demonstrate associativity in the proof of Proposition \ref{cm:prop:monadcondition} may be applied here to show that 
        \[
            \bigcup \left[ \bigcup \Sigma \right]  
                = 
            \bigcup \big[ \, (\Sigma)\ext{\mu}\,  \big]
                \quad \text{for all} \quad
            \Sigma \in \cC \left[ \cflevel{\cC}{n} \left( \moncom{\cC}(S) \right) \right].
    \]
    Altogether, this yields a commuting diagram
    \begin{equation*}
        \begin{tikzcd}
            \cC \left[ \cflevel{\cC}{n} \left( \moncom{\cC}(S) \right) \right]
                \ar[rr, two heads, "\ext{\mu}"]
                \ar[dd, two heads, "\bigcup (-)"']
        &&
            \cC \left[ \moncom{\cC}(S) \right]
                \ar[dd, two heads, "\bigcup (-)"]
        \\
        &&
        \\
            \scK
                \ar[rr, two heads, dashed, "\bigcup (-)"]
                \ar[dr, hook', "i"']
        &&
            \moncom{\cC}(S)
        \\
        &
            \ssc{\scK}
            \ar[ur, two heads, dashed, "\bigcup (-)"']
        &
        \end{tikzcd}
    \end{equation*}
    where the well-definedness of the lower right-hand morphism follows from Lemma \ref{com:lem:downsurj}.
    Since $\ssc{\scK} = \cflevel{\cC}{n+1} \left( \moncom{\cC}(S) \right)$, we are done.
    \item The monad completion map is self-evidently increasing and monotone; and idempotence follows from Lemma \ref{cm:lem:iterateterminate}.
    Consequently, $\moncom{(-)}$ is a closure operator on $\CF$ whose fixed points are precisely the complex monads, and the claim follows.
    \item Since $\cC \leq \moncom{\cC}$, any finite semigroup fixed by $\moncom{\cC}$ must also be a fixed point of $\cC$.
    As for the converse, it is clear that any fixed point of $\cC$ will be fixed by every $\cflevel{\cC}{n}$, and will therefore also be fixed by $\moncom{\cC}$.
\end{enumerate}
\end{proof}

%

%% file: content/mainmatter/catrelmorph.tex
%

\section{The category of relational morphisms}      \label{section:catrelmorph}

\begin{defn}    \label{catrelmorph:defn:rmcat}
The category $\RMCAT$ has the following data.
\begin{enumerate}
    \item Its objects are relational morphisms, and we adopt the convention that a relational morphism $\rho : S \rlm T$ considered as an object of the category $\RMCAT$ will be written as a triple $(S, \rho, T)$.
    \item Its morphisms are of the form
        \[
            (\alpha, \beta) : (S_1, \rho_1, T_1) \arw (S_2, \rho_2, T_2)
        \]
    where $\alpha : S_1 \rightarrow S_2$ and $\beta : T_1 \rightarrow T_2$ are morphisms in $\FSGP$ for which the restricted natural map
    \[
        \alpha \times \beta : \Gamma(S_1, \rho_1, T_1) \arw \Gamma(S_2, \rho_2, T_2)
    \]
    is well-defined in the category $\FSGP$; i.e., such that $(x , y) \in \Gamma(S_1, \rho_1, T_1)$ implies that $(x\alpha , y\beta) \in \Gamma(S_2, \rho_2, T_2)$.
\end{enumerate}
\end{defn}

\subsection{Basic properties}
The relational morphisms 
    \begin{equation*}
        \begin{tikzcd}
            \varnothing \arrow[r] \arrow[d] 
        & 
            \varnothing 
        & 
            \text{and} 
        & 
            \trivsgp \arrow[r] \arrow[d] 
        & 
            \trivsgp 
        \\
            \varnothing 
        & & & 
            \trivsgp 
        &
        \end{tikzcd}
    \end{equation*}
are initial and terminal in $\RMCAT$, respectively.

The category $\RMCAT$ has all finite products, where the product 
    \[
        (S_1 , \rho_1 , T_1) \times (S_2, \rho_2, T_2)
        \;\;=\;\;
        (S_1 \times S_2, \;  \rho_1 \times \rho_2, \; T_1 \times T_2)
    \]
is the relational morphism represented by the factorization
    \begin{equation*}
        \begin{tikzcd}
            \Gamma(S_1 , \rho_1 , T_1) \times \Gamma(S_2, \rho_2, T_2) 
            \arrow[rr, "\rho_{T_1} \times \rho_{T_2}"] 
            \arrow[dd, two heads, "\rho_{S_1} \times \rho_{S_2}"'] 
        && 
            T_1 \times T_2 
        \\&&\\
            S_1 \times S_2 
        &&
        \end{tikzcd}
    \end{equation*}
along with the obvious projections.

In fact, $\RMCAT$ has all finite pullbacks, which are computed via the evident extension of fiber products in $\FSGP$---i.e., the pullback of a cospan
\begin{equation*}
            \begin{tikzcd}
                &&
                    (S_1, \rho_1, T_1) 
                    \ar[dd, "{(\varphi_1, \psi_1)}"]
                \\ && \\
                    (S_2, \rho_2, T_2) 
                    \ar[rr, "{(\varphi_2, \psi_2)}"']
                && 
                    (U, \mu, V)
            \end{tikzcd}
\end{equation*}
in $\RMCAT$ has as its apex the fiber product
    \[
        (S_1, \rho_1, T_1) 
            \times_{(U, \mu, V)}
        (S_2, \rho_2, T_2)
        \;\;=\;\;
        ( S_1 \times_U S_2 , \; \rho_1 \times_\mu \rho_2 , \; T_1 \times_V T_2 )
    \]
whose graph is equal to
    \[
        \Gamma( S_1 \times T_1 , \rho_1 \times \rho_2 , T_1 \times T_2 )
            \cap
        \big[ \left( S_1 \times_U S_2 \right) \times \left( T_1 \times_V T_2 \right) \big].
    \]
Note that the graph of the fiber product in question is isomorphic to the apex of the pullback diagram
\begin{equation*}
    \begin{tikzcd}
                \Gamma(S_1 , \rho_1 , T_1) \times_{\Gamma(U, \mu, V)} \Gamma(S_2, \rho_2, T_2)   
                    \ar[ddd]
                    \ar[rrr]
                    \ar[dddrrr, phantom, "\lrcorner", very near start]
            &&&
                \Gamma(S_1 , \rho_1 , T_1) 
                    \ar[ddd, "\varphi_1 \times \psi_1"]
            \\&&&\\&&&\\
                \Gamma(S_2, \rho_2, T_2) 
                    \ar[rrr, "\varphi_2 \times \psi_2"']
            &&&
                \Gamma(U, \mu, V)
        \end{tikzcd}
\end{equation*}
computed in the category $\FSGP$.

\subsection{}
Associated with $\RMCAT$ is an ensemble cast of functors and natural transformations which are arranged as illustrated by the diagram
\begin{equation*}
    \begin{tikzcd}
        \RMCAT
                \ar[rrrr, bend left=70, "\codRM", ""{name=T, below}]
                    \ar[rrrr, "\Gamma" near start, ""{name=D, above}, ""{name=C, below}]
                \ar[rrrr, bend right=70, "\domRM"', ""{name=B, above}]
        \ar[Rightarrow, from=D, to=T, "\projcod"']
        \ar[Rightarrow, from=C, to=B, "\projdom"]
    &&&&
        \FSGP
    \end{tikzcd}
\end{equation*}
where
\begin{itemize}
    \item the functors $\domRM$, $\Gamma$, and $\codRM$ send a relational morphism to its domain, graph, and codomain, respectively, and their respective actions on morphisms are the obvious ones; and
    \item the components of the natural transformations $\projdom$ and $\projcod$ are the domain and codomain projections, respectively.
\end{itemize}

\begin{defn}
Given $S \in \FSGP$, the category $\RMCAT_S$ has the following data.
\begin{enumerate}
    \item Its objects are relational morphisms of the form $(S, \rho, T)$, i.e., those whose domain is $S$.
    \item Its morphisms are given by $\RMCAT$-morphisms of the form
        \[
            (\id{S} , \psi) : (S, \rho_1, T_1) \arw (S, \rho_2, T_2).
        \]
    We adopt the convention that the $\id{S}$-coordinate of morphisms in $\RMCAT_S$ will be omitted when context allows.
\end{enumerate}
\end{defn}

\subsection{Local products of relational morphisms}
Let $S$ be a finite semigroup.
Given $(S , \rho_1, T_1)$ and $(S, \rho_2, T_2)$ in $\RMCAT_S$, their \textbf{direct sum} is the relational morphism
    \[
        (S, \rho_1 , T_1) \oplus (S , \rho_2, T_2) 
            \, = \, 
        (S , \, \rho_1 \oplus \rho_2, \, T_1 \times T_2)
    \]
whose graph is given by
    \[
        \Gamma(S , \, \rho_1 \oplus \rho_2, \, T_1 \times T_2)
        =
        \setst{(s , (t_1, t_2))}{s \in (t_1)\rho_1^\inv \cap (t_2)\rho_2^\inv}.
    \]
It is easily verified that direct sum is the product in the category $\RMCAT_S$.\footnotemark
\footnotetext{The name "direct sum" for this concept is due to \cite{APRELMOR}.}

Direct sums in $\RMCAT_S$ may be viewed as special pullbacks in $\RMCAT$ as follows.
Given relational morphisms $(S , \rho_1, T_1)$ and $(S, \rho_2, T_2)$ in $\RMCAT_S$, their direct sum (considered as an object of $\RMCAT$) is isomorphic to the fiber product illustrated in the pullback diagram 
\begin{equation*}
    \begin{tikzcd}
                (S, \rho_1, T_1) \times_{(S, \, !_S, \trivsgp)} (S, \rho_2, T_2)   
                    \ar[ddd]
                    \ar[rrr]
                    \ar[dddrrr, phantom, "\lrcorner", very near start]
            &&&
                (S , \rho_1 , T_1) 
                    \ar[ddd, "{(\id{S} ,\, !_{T_1})}"]
            \\&&&\\&&&\\
                (S, \rho_2, T_2) 
                    \ar[rrr, "{(\id{S} ,\, !_{T_2})}"']
            &&&
                (S, \,!_S, \trivsgp)
        \end{tikzcd}
\end{equation*}
where the exclamation marks denote terminal morphisms.

\begin{rmk}
A morphism $\varphi : S_1 \rightarrow S_2$ induces a \textbf{change of base} functor
    \[
        \varphi^\ast : \RMCAT_{S_2} \arw \RMCAT_{S_1}
    \]
whose action, illustrated by
\begin{equation*}
    \begin{tikzcd}
            (S_2, \rho_1, T_1)
                \ar[dd, "\vartheta"', ""{name=A,above}]
    & & & 
                    (S_1, \varphi \rho_1, T_1)
                        \ar[dd, "\vartheta", ""{name=B,below}]
                \ar[from=A, to=B, phantom, "\xmapsto{\qquad \varphi^\ast \qquad}"]
    \\ &&& \\
            (S_2, \rho_2, T_2),        
    & & &           
                    (S_1, \varphi \rho_2, T_2),
    \end{tikzcd}
\end{equation*}
is given by precomposition by $\varphi$ on objects and is "identity" on morphisms.
Note that change of basis preserves direct sums.
The resulting action on morphisms yields a contravariant functor
    \[
        \RMCAT_{(-)} : \FSGP^\op \arw \mathbf{Cat}
    \]
from which the category $\RMCAT$ along with the functor $\domRM$ are recoverable by way of the Grothendieck construction.
\end{rmk}

%

%% file: content/mainmatter/nerves.tex
%

\section{Nerves of relational morphisms}      \label{section:nerves}

\begin{defn}
The \textbf{nerve} of a relational morphism $\rho : S \rlm T$ is the semigroup complex $(S, \Nrv(S, \rho, T))$, where
    \[
        \Nrv(S , \rho , T)
            =
        \ssc{ \setst{ (t)\rho^\inv }{ t \in \Img(\rho) } }.
    \]
In practice, we will refer to both the semigroup complex $(S, \Nrv(S, \rho, T))$ and its face semigroup $\Nrv(S, \rho, T)$ as the ``nerve''.

Moreover, in this situation we say that $\rho : S \rlm T$ \textbf{computes} $\Nrv(S , \rho , T)$.
\end{defn}

\subsection{Alternative descriptions}
The nerve of $\rho : S \rlm T$ consists of sets $X \in \Po(S)$ for which there exists $t \in T$ such that $X \times \{t\}$ is a subset of $\Gamma(S , \rho, T)$.
Another equivalent (and useful) description of the nerve is given by
    \[
        \Nrv(S , \rho , T)
            =
        \setst{ X \in \Po(S) }{ \bigcap_{x \in X} (x)\rho \neq \varnothing }.
    \]

\subsection{The local perspective}
We begin by considering the nerve from a ``local'' or ``pointwise'' perspective; this will prove useful momentarily.
For each finite semigroup $S$ there is a functor
    \[
        \Nrv_S : \RMCAT_S \arw \Com{S},
    \]
where $\Com{S}$ is considered as a category in the standard way, which sends relational morphisms to their nerves and whose action on morphisms is given by
\begin{equation*}
    \begin{tikzcd}
            (S, \rho_1, T_1)
                \ar[dd, "\psi"', ""{name=A,above}]
    & & & 
                    \Nrv(S, \rho_1, T_1)
                        \ar[dd, hook, "\subseteq", ""{name=B,below}]
                \ar[from=A, to=B, phantom, "{\xmapsto{ \quad \Nrv_S \quad }}"]
    \\ &&& \\
            (S, \rho_2, T_2),        
    & & &           
                    \Nrv(S, \rho_2, T_2).
    \end{tikzcd}
\end{equation*}

\begin{prop}    \label{nrv:prop:localnerves}
Let $S$ be a finite semigroup.
\begin{enumerate}
    \item $\Nrv(S , \rho, T) = \sing(S)$ if and only if $\rho : S \sdiv T$ is a division.
    \item The functor $\Nrv_S$ sends direct sums to intersections; that is,
        \[
            \Nrv\big( (S , \rho_1, T_1) \oplus (S , \rho_2, T_2)\big) 
                = 
            \Nrv(S, \rho_1, T_1) \cap \Nrv(S, \rho_2, T_2)
        \]
    for any $(S , \rho_1, T_1), (S, \rho_2, T_2) \in \RMCAT_S$.
\end{enumerate}
\end{prop}

\begin{proof}
Claim (1) is obvious, and claim (2) follows from the observation that
    \begin{align*}
        (t_1 , t_2)[\rho_1 \oplus \rho_2]^\inv
            &= 
        \setst{x \in S}{\textnormal{$(x, t_1) \in \Gamma(S, \rho_1, T_1)$ and $(x, t_2) \in \Gamma(S, \rho_2, T_2)$}}
            \\
            &= 
        (t_1)\rho_1^\inv \cap (t_2)\rho_2^\inv.
    \end{align*}
for any $(t_1, t_2) \in \Img(\rho_1 \oplus \rho_2)$.
\end{proof}

\begin{lem}     \label{nrv:lem:basechange}
If $\varphi : S_1 \rightarrow S_2$ is a morphism, then the diagram
\begin{equation*}
        \begin{tikzcd}
            \RMCAT_{S_2} 
            \ar[d, "\Nrv_{S_2}"'] 
            \ar[rr, "\varphi^\ast"]
        &&
            \RMCAT_{S_1} 
            \ar[d, "\Nrv_{S_1}"]
        \\
            \Com{S_2} 
            \ar[rr, "\plb{\varphi}"]
        && 
            \Com{S_1}
        \end{tikzcd}
\end{equation*}
commutes;
i.e., $\Nrv(S_1 , \varphi \rho, T) = (\Nrv(S_2, \rho, T))\plb{\varphi}$ for any $(S_2, \rho, T) \in \RMCAT_{S_2}$.
\end{lem}

\begin{proof}
Considering the graph
    \[
        \Gamma(S_1, \varphi\rho , T)
            = 
        \setst{(s, t)}{(s\varphi , t) \in \Gamma(S_2, \rho, T)}
    \]
leads one to conclude that $\Nrv(S_1, \varphi\rho , T)$ consists of precisely those sets $X \in \Po(S_1)$ for whom there is some element $t \in T$ such that $(X)\ext{\varphi} \subseteq (t)\rho^\inv$; 
i.e., sets whose image under $\varphi$ is a member of $\Nrv(S_2 , \rho, T)$.
\end{proof}

\subsection{The nerve functor}
Assembling these local pictures yields a functor
    \[
        \comlift{\Nrv} : \RMCAT \arw \COMCAT
    \]
whose action is illustrated by
\begin{equation*}
    \begin{tikzcd}
            (S_1, \rho_1, T_1)
                \ar[dd, "{(\alpha, \beta)}"', ""{name=A,above}]
    & & & 
                    (S_1, \Nrv(S_1, \rho_1, T_1))
                        \ar[dd, "\alpha", ""{name=B,below}]
                \ar[from=A, to=B, phantom, "{\xmapsto{ \quad \comlift{\Nrv} \quad }}"]
    \\ &&& \\
            (S_2, \rho_2, T_2),        
    & & &           
                    (S_2, \Nrv(S_2, \rho_2, T_2)).
    \end{tikzcd}
\end{equation*}
Post-composition of $\comlift{\Nrv}$ by $\comface$ yields a functor $\Nrv$ as in the diagram
\begin{equation*}
            \begin{tikzcd}
            &&
                \COMCAT 
                    \arrow[dd, "\comface"]
            \\&&\\
                \RMCAT 
                    \arrow[uurr, "\comlift{\Nrv}"] 
                    \arrow[rr, "\Nrv"']
            &&
                \FSGP.
            \end{tikzcd}
\end{equation*}

\begin{prop}
The functor $\comlift{\Nrv}$ preserves products; that is, 
    \[
        \Nrv( S_1 \times S_2 , \rho_1 \times \rho_2 , T_1 \times T_2 )
            \;=\;
        \Nrv( S_1 , \rho_1, T_1 ) \otimes \Nrv( S_2, \rho_2, T_2) 
    \]
for any relational morphisms $(S_1 , \rho_1, T_1)$ and $(S_2, \rho_2, T_2)$.
\end{prop}

\begin{proof}
It is easily seen that $(t_1 , t_2)[\rho_1 \times \rho_2]^\inv = (t_1)\rho_1^\inv \times (t_2)\rho_2^\inv$ for any $(t_1, t_2)$ in the image of $\rho_1 \times \rho_2$.
It follows that a face $Z \in \Po(S_1 \times S_2)$ belongs to the nerve of $\rho_1 \times \rho_2$ if and only if 
    \[
        (Z)\ext{\pi_{S_1}} \in \Nrv(S_1, \rho_1, T_1)
            \quad \text{and} \quad
        (Z)\ext{\pi_{S_2}} \in \Nrv(S_2, \rho_2, T_2),
    \]
which is the assertion.
\end{proof}

%

%% file: content/mainmatter/pointlikes.tex
%

\section{Pointlike sets}      \label{section:pointlikes}

\subsection{Codomain specifications}
Given $\pvV \in \PVAR$, the category $\RMCAT^\pvV$ is the full subcategory of $\RMCAT$ spanned by the class of relational morphisms whose codomain belongs to $\pvV$. 
Likewise, $\RMCAT^\pvV_S$ is the analogously defined full subcategory of $\RMCAT_S$ for each $S \in \FSGP$.

\begin{defn}
Let $\pvV \in \PVAR$.
A semigroup complex $(S, \cK)$ is said to be \textbf{$\pvV$-computable} if $\cK = \Nrv(S, \rho, V)$ for some $(S, \rho, V) \in \RMCAT^\pvV$.
\end{defn}

\begin{notation}
Let $\pvV$ be a pseudovariety.
\begin{enumerate}
    \item Let $\CCOMCAT{\pvV}$ denote the full subcategory of $\COMCAT$ spanned by the class of $\pvV$-computable semigroup complexes.
    \item Given $S \in \FSGP$, let $\Com{S}^\pvV$ denote the set of $\pvV$-computable $S$-complexes.
\end{enumerate}
\end{notation}

\subsection{Local closure operators}
Given $\pvV \in \PVAR$ and $S \in \FSGP$, the observation that $\RMCAT^\pvV_S$ is closed under finite direct sums\footnotemark
\footnotetext{Of course, it also implies that $\RMCAT^\pvV$ (and hence $\CCOMCAT{\pvV}$) has all finite products, but this fact is less consequential than its ``local'' counterpart for our purposes.}
implies by way of Proposition \ref{nrv:prop:localnerves} that $\Com{S}^\pvV$ is closed under meets (intersections) in $\Com{S}$.
This yields a closure operator
    \[
        \comclosed{-}{\pvV}_S : \Com{S} \arw \Com{S}
            \quad \text{given by} \quad
        \comclosed{\cK}{\pvV}_S 
            \;=\;
        \bigcap \setst{ \widetilde{\cK} \in \Com{S}^\pvV }{ \cK \subseteq \widetilde{\cK} }
    \]
along with an induced complete lattice structure on $\Com{S}^\pvV$.

\subsection{}
We have finally recovered the ``classical'' notion of pointlikes: 
the bottom element of the lattice $\Com{S}^\pvV$ may be expressed as
    \[
        \comclosed{\sing(S)}{\pvV}_S
            =
        \bigcap \setst{\Nrv(S , \rho , V)}{(S , \rho, V) \in \RMCAT_S^\pvV},
    \]
which is, of course, $\PV(S)$.

\subsection{Computing pointlikes}   \label{existscomputation}
It follows from a straightforward compactness argument that, for any $S \in \FSGP$ and any $\pvV \in \PVAR$, there exists a relational morphism $\rho : S \rlm V$ with $V \in \pvV$ which computes $\PV(S)$; i.e., for which $\Nrv(S, \rho, V) = \PV(S)$.

\begin{lem} \label{pl:lem:closurecompatible}
Let $\varphi : S \rightarrow T$ be a morphism, and let $\cK_S$ and $\cK_T$ be $S$- and $T$-complexes, respectively. 
If $(\cK_S)\psh{\varphi} \subseteq \cK_T$, then $( \comclosed{\cK_S}{\pvV}_S ) \psh{\varphi} \subseteq \comclosed{\cK_T}{\pvV}_T$.
\end{lem}

\begin{proof}
First, note that if $\cK_T'$ is a $\pvV$-computable $T$-complex, then $(\cK_T')\plb{\varphi}$ is $\pvV$-computable as well by Lemma \ref{nrv:lem:basechange}.
Since $\plb{\varphi}$ is an upper adjoint and hence preserves intersections, we have that
    \[
        \left( \comclosed{\cK_T}{\pvV}_T \right) \plb{\varphi}
            \;=\;
        \bigcap \setst{ (\widetilde{\cK_T})\plb{\varphi} }{ \cK_T \subseteq \widetilde{\cK_T} \in \Com{T}^\pvV }.
    \]
By our initial observation, each $(\widetilde{\cK_T})\plb{\varphi}$ in this intersection is $\pvV$-computable. 
Since $\cK_S \subseteq (\widetilde{\cK_T})\plb{\varphi}$ for each of these, it follows that $\comclosed{\cK_S}{\pvV}_S \subseteq ( \comclosed{\cK_T}{\pvV}_T )\plb{\varphi}$, which is equivalent to the desired containment.
\end{proof}

\subsection{From local closure to global closure}
From Lemma \ref{pl:lem:closurecompatible} we obtain a functor
    \[
        \Comp{\pvV} : \COMCAT \arw \CCOMCAT{\pvV}
    \]
whose action is given by
\begin{equation*}
    \begin{tikzcd}
            (S, \cK_S)
                \ar[dd, "\varphi"', ""{name=A,above}]
    & & & 
                    ( S , \comclosed{\cK_S}{\pvV}_S )
                        \ar[dd, "\varphi", ""{name=B,below}]
                \ar[from=A, to=B, phantom, "{\xmapsto{ \quad \Comp{\pvV} \quad }}"]
    \\ &&& \\
            (T, \cK_T),        
    & & &           
                    ( T , \comclosed{\cK_T}{\pvV}_T ).
    \end{tikzcd}
\end{equation*}
The functor $\Comp{\pvV}$ provides a lower adjoint to the evident subcategory inclusion; and, since $\Comp{\pvV}$ is easily seen to be idempotent, the adjunction
\begin{equation*}
            \begin{tikzcd}
            \COMCAT \arrow[rr, shift right=2, "\Comp{\pvV}"', ""{name=B,above}] 
            & &
            \CCOMCAT{\pvV} \arrow[ll, hook', "i"', shift right=2, ""{name=A,below}]
            \ar[from=A, to=B, phantom, "\ddash"]
        \end{tikzcd}
\end{equation*}
exhibits $\CCOMCAT{\pvV}$ as a reflective subcategory of $\COMCAT$.

\subsection{Pointlike factorization}
For each $\pvV \in \PVAR$, the $\pvV$-pointlikes functor $\PV$ may be defined as the composition indicated by the diagram
\begin{equation*}
        \begin{tikzcd}
            \COMCAT  
            \ar[rr, "\Comp{\pvV}"]
        &&
            \CCOMCAT{\pvV} 
            \ar[dd, "\comface"]
        \\&&\\
            \FSGP
            \ar[uu, "\comlift{\sing}"]
            \ar[rr, dashed, "\PV"]
        && 
            \FSGP.
        \end{tikzcd}
\end{equation*}
Considering Lemma \ref{pl:lem:closurecompatible} in the context of Proposition \ref{com:prop:monoepiequiv} immediately shows that $\PV$ preserves regular epimorphisms, which in turn implies that $\PV$ is a complex functor (since the other two conditions are obviously satisfied).
It is similarly easy to see that $\Fix(\PV) = \pvV$.
In order to show that $\PV$ is the maximum complex functor fixing $\pvV$, we require the following lemma.

\begin{lem} \label{lem:fixedpointnerve}
Let $\cC$ be a complex functor.
If $\rho : S \rlm T$ is a relational morphism and $T$ is a fixed point of $\cC$, then $\cC(S) \subseteq \Nrv(S , \rho, T)$.
\end{lem}

\begin{proof}
Let $\ext{\rho} : \cC(S) \rlm \cC(T)$ denote the relational morphism obtained by applying $\cC$ to the canonical factorization of $\rho$.
\begin{equation*}
    \begin{tikzcd}
    && 
        \Gamma(\cC(S) , \ext{\rho} , \cC(T))
            \ar[dll, two heads, bend right=15]
            \ar[drr, bend left=15]
    &&
    \\
        \cC(S)
    &&
        \cC(\Gamma(S , \rho, T))
            \ar[ll, two heads]
            \ar[rr]
            \ar[u, two heads]
    &&
        \cC(T)
    \\
        S
            \ar[u, hook]
    &&
        \Gamma(S , \rho, T)
            \ar[ll, two heads]
            \ar[rr]
            \ar[u, hook]
    &&
        T
            \ar[u, hook]
    \end{tikzcd}
\end{equation*}
Now, if $Z \in \cC(\Gamma(S , \rho, T))$ then the set
    \[
        (Z)\ext{\rho_T} \:=\: \setst{t}{(s, t) \in Z}
    \]
is a singleton ex hypothesi, and so the graph of $\ext{\rho}$ consists of pairs $(X, \{t\})$ for which $X \times \{t\}$ is a subset of the graph of $\rho$.
Since every such $X$ is a member of $\Nrv(S, \rho, T)$, surjectivity of $\ext{\rho_S}$ yields the lemma.
\end{proof}

\begin{thm}     \label{MAINTHEOREM}
\hfill
\begin{enumerate}
    \item If $\pvV \in \PVAR$, then $\PV$ is maximal amongst complex functors fixing $\pvV$; that is,
        \[
            \PV = \bigvee \setst{ \cC \in \CF }{ \pvV \subseteq \Fix ( \cC ) }.
        \]
    Consequently, the map
        \[
            \PLF : \PVAR^\op \arw \CF 
                \quad \text{given by} \quad
            \pvV \longmapsto \PV
        \]
    is upper adjoint to $\Fix$.
    Additionally, $\PV$ is a complex monad.
    \item Altogether, there is a commuting triangle 
        \begin{equation*}
        \begin{tikzcd}
            \CF
            \ar[rrrr, shift right=1, bend left=15, ""{name=B,above}, two heads, "\moncom{(-)}"']
            \ar[dddrr, shift left=1, bend right=15, ""{name=C}, two heads, "\Fix"]
            &       &       &       &
            \CM
            \ar[llll, shift right=3, bend right=15, ""{name=A,below}, hook', "i"']
            \ar[dddll, shift right=1.5, bend left=15, two heads, ""{name=E}, "\Fix"']
            \ar[from=A, to=B, phantom, "\scriptstyle{\ddash}"]
        \\
            &       &       &       &
        \\
            &       &       &       &
        \\
            &       &
            \PVAR^\op
            \ar[uuull, shift left=3, bend left=15, ""{name=D}, hook, "\PLF"]
            \ar[uuurr, shift right=3, bend right=15, ""{name=F}, hook', "\PLF"']
            &       &
        \ar[from=C, to=D, phantom, "\scriptstyle{\ddash}" rotate=135]
        \ar[from=E, to=F, phantom, "\scriptstyle{\ddash}" rotate=-135]
        \end{tikzcd}
        \end{equation*}
    of Galois connections between complete lattices.
\end{enumerate}
\end{thm}

\begin{proof}
\hfill
\begin{enumerate}
    \item Let $\cC \in \CF$ with $\pvV \subseteq \Fix(\cC)$.
    As noted in \ref{existscomputation}, for any $S \in \FSGP$ we are guaranteed the existence of some $\rho : S \rlm V$ with $V \in \pvV$ for which $\Nrv(S, \rho, V) = \PV(S)$.
    Applying Lemma \ref{lem:fixedpointnerve} shows that $\cC(S) \subseteq \PV(S)$ and hence that $\cC \leq \PV$.
    This establishes the claimed maximality of $\PV$, and the claim regarding $\PLF$ follows.
    
    The claim that $\PV$ is a complex monad follows by maximality from claims (2) and (3) of Proposition \ref{SECONDPUNCHLINE}, which state that the monad completion map $\moncom{(-)}$ is a closure operator which preserves fixed points.
    \item This is a summary obtained by considering claim (1) alongside Propositions \ref{FIRSTPUNCHLINE} and \ref{SECONDPUNCHLINE}.
\end{enumerate}
\end{proof}

%

%% file: content/mainmatter/moduli.tex
%

\section{Moduli}      \label{section:moduli}

\begin{defn}    \label{mod:defn:moduli}
A \textbf{modulus} is a rule $\Lambda$, written as 
    \[
        \Lambda = \modrl{S}{\Lambda_S},
    \]
which assigns to each finite semigroup $S$ a (possibly empty) set $\Lambda_S \subseteq \Po(S)$, and which satisfies the following axioms.
\begin{enumerate}
    \item If $\varphi : S \rightarrow T$ is a morphism, then for any $X \in \Lambda_S$ there exists some $\widetilde{X} \in \Lambda_T$ such that $(X)\ext{\varphi} \subseteq \widetilde{X}$.
    \item If $\varphi : S \surj T$ is a regular epimorphism, then for any $Y \in \Lambda_T$ there exists some $\widetilde{Y} \in \Lambda_S$ such that $(\widetilde{Y})\ext{\varphi} = Y$.
\end{enumerate}
\end{defn}

\begin{example} \label{moduliexamples}
Commonplace moduli include
\begin{enumerate}
    \item the \textit{subgroup} modulus
        \[
            \mathsf{Grp} = \modrl{S}{ \setst{G}{ \textnormal{$G$ is a subgroup of $S$} } };
        \]
    \item the three \textit{Green's} moduli\footnotemark
    \footnotetext{Note that the analogous definition for $\GH$ does \textit{not} yield a modulus (see \cite[Chapter 7]{ARBIB} or \cite{CONDITIONALEQ}).}    
        \[
        \mathsf{RCl} = \modrl{S}{S / \GR},
            \quad
        \mathsf{LCl} = \modrl{S}{S / \GL},
            \quad \text{and} \quad
        \mathsf{JCl} = \modrl{S}{S / \GJ};
        \]
    \item the three ``principal'' moduli, including the \textit{principal right ideals} modulus
        \[
            \mathsf{PrinR}
                =
            \modrl{S}{ \setst{ x \cdot S^I }{x \in S} },
        \]
    the \textit{principal left ideals} modulus
        \[
            \mathsf{PrinL}
                =
            \modrl{S}{ \setst{S^I \cdot x}{x \in S} },
        \]
    and the \textit{principal (two-sided) ideals} modulus
        \[
            \mathsf{PrinJ}
                =
            \modrl{S}{ \setst{S^I \cdot x \cdot S^I}{x \in S} };
        \]
    \item given a positive integer $k$, the \textit{$k$-length products} modulus
        \[
            \mathsf{Prod}^k 
                =
            \modrl{S}{ \big\{ \setst{ x_1 x_2 \cdots x_k }{ x_i \in S} \big\} };
        \]
    \item again for a positive integer $k$, the \textit{$k$-length suffix} modulus
        \[
            \mathsf{Suffix}^k 
                =
            \modrl{S}{ \setst{ S^I \cdot (x_1 x_2 \cdots x_k) }{ x_i \in S}  }
        \]
    and the \textit{$k$-length prefix} modulus
        \[
            \mathsf{Prefix}^k 
                =
            \modrl{S}{ \setst{ (x_1 x_2 \cdots x_k) \cdot S^I }{ x_i \in S}  };
        \]
    \item the \textit{set of idempotents} modulus 
        \[
            \mathsf{E} = \modrl{S}{\{E(S)\}};
        \]
    \item and the \textit{set of regular elements} modulus
        \[
            \mathsf{Reg} = \modrl{S}{\{\Reg(S)\}}.
        \]
\end{enumerate}
\end{example}

\begin{rmk}
Moduli are similar to the ``implicit relations'' considered in \cite{CONDITIONALEQ}.
\end{rmk}

\begin{notation}
Let $\MD$ denote the collection of moduli.
\end{notation}

\subsection{Order theoretic aspects}
A modulus $\Lambda_1$ is said to \textbf{refine} another modulus $\Lambda_2$ if, for all finite semigroups $S$, any set $X_1 \in \Lambda_{1 , S}$ is a subset of some set $X_2 \in \Lambda_{2 , S}$.
Refinement is a preorder on $\MD$, and the moduli $\modrl{S}{\Po(S)}$ and $\modrl{S}{\varnothing}$ are respectively maximal and minimal with respect to refinement.

Joins exist in $\MD$ with respect to refinement and are given by
    \[
        \Lambda_1 \vee \Lambda_2
        \:=\: \modrl{S}{\Lambda_{1 , S} \cup \Lambda_{2 , S}}
    \]
for any moduli $\Lambda_1$ and $\Lambda_2$.

\begin{defn}
The set of \textbf{points} of a modulus $\Lambda$ is 
    \[
        \points{\Lambda}    
            \;=\;
        \setst{S \in \FSGP}{\Lambda_S \subseteq \sing(S)}.
    \]
\end{defn}

\begin{prop}
If $\Lambda$ is a modulus, then $\points{\Lambda}$ is a pseudovariety.
\end{prop}

\begin{proof}
First, let $\varphi : S \inj T$ be a monomorphism with $T \in \points{\Lambda}$.
Then for any $X \in \Lambda_S$, its image $(X)\ext{\varphi}$ must be a subset of some singleton in $\Lambda_T$, which, since $\varphi$ is injective, implies that $X$ must be a singleton as well; from which it follows that $S \in \points{\Lambda}$ as well.

Next, let $\varphi : S \surj T$ be a regular epimorphism with $S \in \points{\Lambda}$.
Then for any $Y \in \Lambda_T$ there exists some singleton in $\Lambda_S$ whose image under $\varphi$ is $Y$; hence $Y$ must also be a singleton, and consequently $T \in \points{\Lambda}$ as well.

Finally, consider $S, T \in \points{\Lambda}$.
Given $Z \in \Lambda_{S \times T}$, there exist $Z_S \in \Lambda_S$ and $Z_T \in \Lambda_T$ such that $(Z)\ext{\pi_S} \subseteq Z_S$ and $(Z)\ext{\pi_T} \subseteq Z_T$.
But, since both $Z_S$ and $Z_T$ must be singletons, $Z$ must be a singleton as well.
Thus $S \times T \in \points{\Lambda}$ as well, and we are done.
\end{proof}

\begin{example}
Consider again the moduli in Example \ref{moduliexamples}.
\begin{enumerate}
    \item Clearly $\points{ \mathsf{Grp} } = \pvA$.
    \item The points of the three Green's moduli are given by
        \[
        \points{\mathsf{RCl}} = \pvR,
            \quad
        \points{\mathsf{LCl}} = \pvL,
            \quad \text{and} \quad
        \points{\mathsf{JCl}} = \pvJ.
        \]
    \item The points of the principal moduli are given by
        \[
        \points{\mathsf{PrinR}} = \mathbf{LZ},
            \quad
        \points{\mathsf{PrinL}} = \mathbf{RZ},
            \quad \text{and} \quad
        \points{\mathsf{PrinJ}} = \pvtriv.
        \]
    where $\mathbf{LZ}$ and $\mathbf{RZ}$ are the pseudovarieties of left- and right-zero semigroups, respectively.
    \item The points of the $k$-length product modulus are given by
        \[
            \points{\mathsf{Prod}^k}
                =
            \pvN_k
        \]
    where $\pvN_k$ is the pseudovariety of $k$-nilpotent semigroups.
    \item The points of the $k$-length suffix and prefix moduli are given by
        \[
            \mathsf{Suffix}^k = \pvD_k
                \quad \text{and} \quad
            \mathsf{Prefix}^k = \pvK_k,
        \]
    respectively, where $\pvD_k$ is the \textit{level $k$ delay} pseudovariety whose members satisfy the equation
        \[
            y (x_1 x_2 \cdots x_k)  = (x_1 x_2 \cdots x_k),
        \]
    and $\pvK_k$ is the \textit{level $k$ reverse delay} pseudovariety whose members satisfy the equation
        \[
            (x_1 x_2 \cdots x_n) y  = (x_1 x_2 \cdots x_n).
        \]
    \item Since a finite semigroup belongs to $\points{ \mathsf{E} }$ if and only if it has a unique idempotent, one has that      
    \[
        \points{ \mathsf{E} } = \pvG \malcev \pvN.
    \]
    \item Similarly, a finite semigroup has a unique regular element if and only if it is nilpotent, and so
    \[
        \points{ \mathsf{Reg} } = \pvN.
    \]
\end{enumerate}
\end{example}

\begin{defn}
A modulus $\Lambda$ induces a complex functor $\modcf{\Lambda}$ given by
    \[
        \modcf{\Lambda}(S)
            \;=\;
        \COMGEN{ \Lambda_S }{S}
    \]
at every finite semigroup $S$.
\end{defn}

\begin{prop} \label{mod:prop:construct}
\hfill
\begin{enumerate}
    \item The map 
    \[
        \sfC : \MD \arw \CF
            \quad \text{given by} \quad
        \Lambda \longmapsto \modcf{\Lambda}(S)
    \]
    is monotone and join-preserving.
    \item $\Fix(\modcf{\Lambda}) = \points{\Lambda}$ for any modulus $\Lambda$.
\end{enumerate}
\end{prop}

\begin{proof}
Straightforward.
\end{proof}

\begin{notation}
Given $\Lambda \in \MD$, let $\modcm{\Lambda}$ denote the monad completion of $\modcf{\Lambda}$.
\end{notation}

\begin{thm} \label{MODULARCONSTRUCTIONTHEOREM}
Let $\pvV$ be a pseudovariety and let $\Lambda$ be a modulus. Then
    \[
        \pvV \subseteq \points{\Lambda}
        \quad \Longleftrightarrow \quad
        \modcm{\Lambda} \leq \PV.
    \] 
\end{thm}

\begin{proof}
Considering claim (2) of Proposition \ref{mod:prop:construct} in the context of Theorem \ref{MAINTHEOREM} yields the theorem.
\end{proof}

%

%% file: content/mainmatter/effectiveness.tex
%

\section{Effectiveness of moduli}      \label{section:effectiveness}

\subsection{}
Theorem \ref{MODULARCONSTRUCTIONTHEOREM} formalizes a ubiquitous theme in the study of pointlike sets wherein lower bounds for pointlikes are constructed by iterative unioning of distinguished subsets.

More precisely, it reduces the problem of constructing lower bounds for $\pvV$-pointlikes to the problem of defining a modulus $\Lambda$ for which $\Lambda_V$ is at most singletons for any $\pvV$-semigroup $V$.
Once such a modulus $\Lambda$ is obtained, it remains to show that $\scC_\Lambda$ is also an upper bound for $\PV$---a task which is, in general, considerably more difficult. 

\begin{defn}
A modulus $\Lambda$ is said to be \textbf{effective} with respect to a pseudovariety $\pvV$---or \textit{$\pvV$-effective} for short---if $\scC_\Lambda = \PV$.
\end{defn}

\begin{example}[Aperiodics; see \cite{PL-APERIODIC-KARSTEN, PRODEXP, QTHEORY}] \label{plexample:aperiodic}
The modulus
    \[
        \mathsf{Grp} = \modrl{S}{ \setst{G}{ \textnormal{$G$ is a subgroup of $S$} } }
    \]
is effective with respect to the pseudovariety $\pvA$ of aperiodic semigroups.
An alternative $\pvA$-effective modulus appearing in the literature is
    \[
        \mathsf{CycGrp} = \modrl{S}{ \setst{ \langle g \rangle }{ \textnormal{$g$ is a group element of $S$} } }.
    \]
\end{example}

\begin{example}[$\GR$- and $\GL$-trivial; see \cite{HYP-R-ALM, PL-RJ-ALM}] \label{plexample:RandLtrivial}
The moduli
    \[
        \mathsf{RCl} = \modrl{S}{S/\GR}
            \quad \text{and} \quad
        \mathsf{LCl} = \modrl{S}{S/\GL}
    \]
are effective with respect to the pseudovarieties $\pvR$ and $\pvL$ of $\GR$- and $\GL$-trivial semigroups, respectively.
\end{example}

\begin{rmk} \label{existenceofimposters} 
Not every modulus is effective with respect to its pseudovariety of points; and, consequently, the Galois connection between $\CM$ and $\PVAR^\op$ is \textit{not} an equivalence. 
To see this, consider the modulus
    \[
        \mathsf{JCl} = 
        \modrl{S}{S / \GJ}.
    \]
Clearly $\points{\mathsf{JCl}} = \pvJ$, and hence $\modcm{\mathsf{JCl}} \leq \PL{\pvJ}$ by Theorem \ref{MODULARCONSTRUCTIONTHEOREM}.

However, there is a semigroup constructed in \cite[Subsection~6.1]{PL-RJ-ALM} which---when interpreted in the language of this paper---provides an example of a finite semigroup $S$ for whom the containment $\modcm{\mathsf{JCl}}(S) \subset \PL{\pvJ}(S)$ is \textit{strict}; and so the inequality $\modcm{\mathsf{JCl}} < \PL{\pvJ}$ is strict as well.
\end{rmk}

%

%% file: content/mainmatter/transfer.tex
%

\section{A framework for transfer results}      \label{section:transfer}

\subsection{}
Our approach to transfer results begins with two diagrams 
    \begin{equation*}
        \begin{tikzcd}
            \CF 
                \ar[rr, two heads, "\Fix"] 
                \ar[dd, dashed] 
        && 
            \PVAR
                \ar[dd, dashed]
        && 
            \CF  
                \ar[dd, dashed] 
        && 
            \PVAR
                \ar[ll, hook', "\PLF"']
                \ar[dd, dashed]
        \\&&& \text{and} &&&\\
            \CF
                \ar[rr, two heads, "\Fix"]
        && 
            \PVAR
        && 
            \CF
        &&
            \PVAR
                \ar[ll, hook', "\PLF"']
        \end{tikzcd}
    \end{equation*}
corresponding to the lower and upper adjoints of the main Galois connection.
The program, broadly speaking, is concerned with pairs of operators 
    \[
        \lambda : \CF \arw \CF
            \quad \text{and} \quad
        \alpha : \PVAR \arw \PVAR
    \]
whose placement along the dashed arrows makes one or both diagrams commute.

\begin{notation}
In this section, the application of ``generic'' operators on the lattices $\CF$ and $\PVAR$ will be written with superscripts; e.g., for $\lambda$ and $\alpha$ as above, write
    \[
        \lambda : \cC \longmapsto \cC^\lambda
            \quad \text{and} \quad
        \alpha : \pvV \longmapsto \pvV^\alpha
    \]
for any $\cC \in \CF$ and any $\pvV \in \PVAR$.
\end{notation}

\subsection{Classes of operators}
We begin by considering what sort of operators we are interested in. 
In general, the relevant operators on $\PVAR$ are those which are continuous. 
This has long been established---however, there is no such precedent with respect to operators on $\CF$. 
It would be worthwhile to study various classes of and conditions on operators on $\CF$, but we will not do so in any depth here. 
For our purposes, we will consider a very general class of operators on CF given by a single (and natural) condition.

\begin{defn}
A monotone operator $\lambda : \CF \rightarrow \CF$ is said to be \textbf{algebraic} if
    \[
        \Fix(\cC_1) = \Fix(\cC_2)
            \quad \Longrightarrow \quad
        \Fix(\cC_1^\lambda) = \Fix(\cC_2^\lambda)
    \]
for all $\cC_1 , \cC_2 \in \CF$.
\end{defn}

\begin{notation}
Let $\AlgCF$ denote the collection of algebraic operators on $\CF$.
\end{notation}

\begin{defn}    \label{tr:defn:conditions}
Let $(\lambda, \alpha) \in \AlgCF \times \CntPV$.
\begin{enumerate}
    \item The pair $(\lambda, \alpha)$ satisfies the \textbf{fixed point transfer condition} (FPTC) if
    \begin{equation*}
        \begin{tikzcd}
            \CF 
                \ar[d, "\lambda"'] 
                \ar[rr, two heads, "\Fix"]
        &&
            \PVAR 
                \ar[d, "\alpha"]
        \\
            \CF 
                \ar[rr, two heads, "\Fix"']
        && 
            \PVAR
        \end{tikzcd}
    \end{equation*}
    is commutative.
    \item Likewise, if the diagram
    \begin{equation*}
        \begin{tikzcd}
            \CF 
                \ar[d, "\lambda"'] 
        &&
            \PVAR
                \ar[ll, hook', "\PLF"']
                \ar[d, "\alpha"]
        \\
            \CF 
        && 
            \PVAR
                \ar[ll, hook', "\PLF"']
        \end{tikzcd}
    \end{equation*}
    commutes, then $(\lambda, \alpha)$ satisfies the \textbf{pointlike transfer condition} (PLTC).
\end{enumerate}
\end{defn}

\begin{notation}
Let $\FPTC$ and $\PLTC$ denote the collections of pairs satisfying the fixed point and pointlike transfer conditions, respectively.
\end{notation}

\subsection{}
It is straightforwardly seen that $\PLTC \subseteq \FPTC$, and it is similarly evident that both $\PLTC$ and $\FPTC$ are closed under coordinatewise composition.

\subsection{}
Since, in practice, we often have a particular operator $\alpha \in \CntPV$ in mind, it is convenient to say that $\lambda \in \AlgCF$ \textit{satisfies FPTC with respect to $\alpha$} to mean that the pair $(\lambda, \alpha)$ satisfies FPTC.
We also adopt the analogous convention in the case of PLTC.

\subsection{}
Proving satisfaction of the pointlike transfer condition is likely to be extremely difficult in most cases.
Hence, in the interest of setting realistic goals, it is natural to consider ``restricted'' satisfaction of this condition.
To this end, given a class $\scV \subseteq \PVAR$, a pair $(\lambda, \alpha) \in \AlgCF \times \CntPV$ for which the diagram
\begin{equation*}
            \begin{tikzcd}
                    \CF 
                    \ar[d, "\lambda"'] 
                &&
                    \scV 
                    \ar[d, "\alpha"]
                    \ar[ll, hook', "\PLF"']
                \\
                    \CF 
                && 
                    \PVAR
                    \ar[ll, hook', "\PLF"]
            \end{tikzcd}
\end{equation*}
commutes satisfies the \textbf{local pointlike transfer condition on $\scV$}.

%

%% file: content/mainmatter/reversal.tex
%

\section{Example: pointlike transfer along reversal}      \label{section:reversal}

\subsection{}
The \textbf{reversal} of $S \in \FSGP$ is the semigroup $S^\rev$ whose underlying set is that of $S$ but whose operation is \textit{reversed} in the sense that, for any $x, y \in S^\rev$,
    \[
        x \cdot^\rev y = yx,
    \]
where $\cdot^\rev$ denotes the product in $S^\rev$.

Reversal is an endofunctor in the evident way, and moreover this functor may be extended to an operator 
    \[
        \rev : \PVAR \arw \PVAR    
            \quad \text{given by} \quad
        \pvV^\rev = \setst{V^\rev}{V \in \pvV}
    \]
for every pseudovariety $\pvV$.
Both the functor and the operator personae are involutions on their respective domains.

\subsection{Conjugation by reversal}
Consider the operator
    \[
        \rightleftarrows \; : \CF \arw \CF
            \quad \text{given by} \quad
        \cC^\rightleftarrows(S) = \big( \cC(S^\rev) \big)^\rev
    \]
for every $\cC \in \CF$ and every $S \in \FSGP$.
In other words, the operator $\rightleftarrows$ acts via conjugation by $\rev$ as illustrated by the diagram
\begin{equation*}
    \begin{tikzcd}
        \FSGP 
            \ar[dd, dashed, "\cC^\rightleftarrows"'] 
            \ar[rr, "\rev"]
    &&
        \FSGP 
            \ar[dd, "\cC"]
    \\&&\\
        \FSGP 
    && 
        \FSGP.
            \ar[ll, "\rev"]
    \end{tikzcd}
\end{equation*}
It is obvious that $\rightleftarrows$ is algebraic.

\begin{prop} \label{rev:prop:reversalpltc}
The pair $(\rightleftarrows , \rev)$ satisfies the pointlike transfer condition.
That is, 
    \[
        \PL{\pvV^\rev}(S) = \big( \PV(S^\rev) \big)^\rev
    \]
for any pseudovariety $\pvV$ and any finite semigroup $S$.
\end{prop}

\begin{proof}
Extend reversal to relational morphisms by defining $(S , \rho , T)^\rev$ to be the relational morphism $(S^\rev , \rho^\rev , T^\rev)$ whose graph is given by
    \[
        \Gamma \big( S^\rev , \rho^\rev , T^\rev \big) 
            \:\cong\: 
        \big( \Gamma(S , \rho , T) \big)^\rev,
    \]
i.e., the graph of $\rho^\rev$ is isomorphic to the reversal of the graph of $\rho$.
It is easily verified that
    \[
        \Nrv \big( S^\rev , \rho^\rev , T^\rev \big) 
            \:=\: 
        \big( \Nrv(S , \rho, T) \big)^\rev.
    \]

Now, if $\rho : S^\rev \rlm V$ computes $\PV(S^\rev)$, then $\big( S, \rho^\rev , V^\rev \big)$ is a relational morphism belonging to $\RMCAT_S^{\pvV^\rev}$; hence
    \[
        \PL{\pvV^\rev}(S) 
        \;\;\subseteq\;\; \Nrv \big( S, \rho^\rev , V^\rev \big) 
        \;\;=\;\; \big( \PV(S^\rev) \big)^\rev.
    \]
The symmetric argument concerning a relational morphism $\rho : S \rlm V^\rev$ computing $\PL{\pvV^\rev}(S)$ yields the proposition.
\end{proof}

%

%% file: content/mainmatter/contextspec.tex
%

\section{Context specifiers}      \label{section:contextspec}

\begin{defn}    \label{cts:defn:contextspec}
A \textbf{context specifier} is a rule $\scO$, written as
    \[
        \scO \;=\; \modrl{S}{\scO(S)},
    \]
which assigns to each finite semigroup $S$ a (possibly empty) set $\scO(S)$ of subsemigroups of $S$ in a manner which satisfies the following axioms.
\begin{enumerate}
    \item For all morphisms $\varphi : S \rightarrow T$ and all $U_S \in \scO(S)$ there exists some $U_T \in \scO(T)$ for which the restriction of $\varphi$ to $U_S$ factors through the inclusion of $U_T$ into $T$.
    \begin{equation*}
        \begin{tikzcd}
            U_S 
                \ar[d, hook] 
                \ar[rr, "\varphi \mid_{U_S}", dashed]
        &&
            U_T 
                \ar[d, hook, dashed, "\exists"]
        \\
            S
                \ar[rr, "\varphi"]
        && 
            T
        \end{tikzcd}
    \end{equation*}
    \item For all regular epimorphisms $\varphi : S \surj T$ and all $U_T \in \scO(T)$ there exists some $U_S \in \scO(S)$ for which the restriction of $\varphi$ to $U_S$ is a regular epimorphism with image $U_T$.
    \begin{equation*}
        \begin{tikzcd}
            U_S 
                \ar[d, hook, dashed, "\exists"'] 
                \ar[rr, two heads, "\varphi \mid_{U_S}", dashed]
        &&
            U_T 
                \ar[d, hook]  
        \\
            S
                \ar[rr, two heads, "\varphi"]
        && 
            T
            \end{tikzcd}
    \end{equation*}
\end{enumerate}
\end{defn}

\begin{notation}
The collection of context specifiers is denoted by $\CTS$.
\end{notation}

\begin{example}
Familiar context specifiers include
\begin{enumerate}
\item the \textit{subgroup} context
    \[
        \mathsf{Grp} = \modrl{S}{ \setst{G}{ \textnormal{$G$ is a subgroup of $S$} } }
    \]
and the \textit{cyclic subgroup} context
    \[
        \mathsf{CycGrp} = \modrl{S}{ \setst{ \langle g \rangle }{ \textnormal{$g$ is a group element of $S$} } }
    \]
as considered in Example \ref{plexample:aperiodic};
\item the \textit{local} context 
    \[
        \mathsf{Loc} = \modrl{S}{ \setst{eSe}{ e \in E(S) } };
    \]
\item given two pseudovarieties $\pvV$ and $\pvW$, the context 
    \[
        \mathsf{\pvV-Like}_\pvW 
            =
        \modrl{S}{ \setst{ T \leq S }{ \text{$T$ is $\pvV$-like with respect to $\pvW$} } },
    \]
where a subsemigroup $T$ of $S$ is said to be \textit{$\pvV$-like with respect to $\pvW$} if for every relational morphism $\rho : S \rlm W$ with $W \in \pvW$, there exists some $\pvV$-subsemigroup $V$ of $W$ for which $T \subseteq (V)\rho^\inv$;
\item given a pseudovariety $\pvW$, the \textit{$W$-idempotent pointlike} context
        \[
            \mathsf{EP}_\pvW = \mathsf{\pvtriv-Like}_\pvW,
        \]
which assigns $S \in \FSGP$ to its set of $\pvW$-idempotent pointlike subsemigroups, which are subsemigroups $T$ of $S$ such that whenever $\rho : S \rlm W$ is a relational morphism with $W \in \pvW$, there exists an idempotent $e \in W$ for which $T \subseteq (e)\rho^\inv$;
\item and, given a pseudovariety $\pvH$ of groups, the \textit{$\pvH$-kernel} context
    \[
        \mathsf{Ker}_\pvH = \modrl{S}{ \{ K_\pvH(S) \} },
    \]
where $K_\pvH(S)$, called the \textit{$\pvH$-kernel} of $S$, is the unique maximal $\pvH$-idempotent pointlike subsemigroup of $S$.
\end{enumerate}
\end{example}

\subsection{Generation by moduli}
Context specifiers are a special case of moduli, and the inclusion of $\CTS$ into $\MD$ is a section of the map
    \[
        \CTSGEN{-} : \MD \longsurj \CTS 
        \quad \text{given by} \quad
        \CTSGEN{\Lambda} = \modrl{S}{ \setst{ \langle X \rangle }{ X \in \Lambda_S } },
    \]
which is a monotone retraction of $\MD$ onto $\CTS$.

\begin{example}
Commonplace context specifiers arising from $\CTSGEN{-}$ include
\begin{enumerate}
\item the \textit{idempotent generated} context
    \[
        \mathsf{EGen} = \CTSGEN{\mathsf{E}} = \modrl{S}{ \{ \langle E(S) \rangle \} };
    \]
\item and the \textit{regular generated} context
    \[
        \mathsf{RegGen} = \CTSGEN{\mathsf{Reg}} = \modrl{S}{ \{ \langle \Reg(S) \rangle \} }.
    \]
\end{enumerate}
\end{example}

\begin{defn} \label{cts:defn:contextgalois}
A context specifier $\scO$ induces two operators
\begin{align*}
            \CTSinPV{\scO}{-} : \PVAR \arw \PVAR
            & \quad\quad \text{given by} \quad\quad
            \CTSinPV{\scO}{\pvV} 
                \;=\;
            \setst{ S \in FSGP }{ \scO(S) \subseteq \pvV }
\\
& \quad\quad\;\;\; \text{and} 
\\
            \CTSsumPV{\scO}{-} : \PVAR \arw \PVAR
            & \quad\quad \text{given by} \quad\quad
            \CTSsumPV{\scO}{\pvV}
                \;=\;
            \PVGEN{ \bigcup \setst{\scO(V)}{V \in \pvV} }.
\end{align*}
\end{defn}

\begin{prop}
Given a context specifier $\scO$, the operators given in Definition \ref{cts:defn:contextgalois} are the constituents of a Galois connection
\begin{equation*}
            \begin{tikzcd}
            \PVAR \arrow[rrr, "\CTSsumPV{\scO}{-}"', bend right=12, ""{name=B,above}] 
            & & &
            \PVAR. \arrow[lll, "\CTSinPV{\scO}{-}"', bend right=12, ""{name=A,below}]
            \ar[from=A, to=B, phantom, "\ddash"]
        \end{tikzcd}
\end{equation*}
Moreover, $\CTSinPV{\scO}{-}$ is continuous and increasing, and $\CTSsumPV{\scO}{-}$ is decreasing.
\end{prop}

\begin{proof}
It is clear that both $\CTSinPV{\scO}{-}$ and $\CTSsumPV{\scO}{-}$ are monotone, and it is similarly clear that the former is increasing and the latter is decreasing.

To see that these maps are adjoint, let $\pvV$ and $\pvW$ be pseudovarieties.
Then $\CTSsumPV{\scO}{\pvV} \subseteq \pvW$ if and only if $\scO(V) \subseteq \pvW$ for every $V \in \pvV$, which is the case if and only if $\pvV \subseteq \CTSinPV{\scO}{\pvW}$.

Finally, to see that $\CTSinPV{\scO}{-}$ is continuous, let $\scV = \{ \pvV_\alpha \}_{\alpha \in I}$ be a directed set of pseudovarieties.
Note that monotonicity implies that the image of $\scV$ under $\CTSinPV{\scO}{-}$ is directed as well, and hence that
    \[
        \bigvee_{\alpha \in I} \CTSinPV{\scO}{ \pvV_\alpha }
            \;=\;
        \bigcup_{\alpha \in I} \CTSinPV{\scO}{ \pvV_\alpha }.
    \]
Now, a finite semigroup $S$ belongs to $\CTSinPV{\scO}{ \bigcup \scV }$ if and only if $\scO(S) \subseteq \bigcup \scV$, which is the case if and only if $S$ belongs to $\bigcup_{\alpha \in I} \CTSinPV{\scO}{ \pvV_\alpha }$.
\end{proof}

\begin{example}
Familiar context specifiers induce familiar operators:
\begin{enumerate}
\item the operators induced by $\mathsf{Grp}$ are  
    \[
        \CTSinPV{\mathsf{Grp}}{-} = \mathbb{G}[-]
            \qquad \text{and} \qquad
        \CTSsumPV{\mathsf{Grp}}{-} = \pvG \cap (-);
    \]
\item the operators induced by $\mathsf{Loc}$ are 
    \[
        \CTSinPV{\mathsf{Loc}}{-} \;=\; \mathbb{L}[-]
            \qquad \text{and} \qquad
        \CTSsumPV{\mathsf{Loc}}{-} \;=\; \PVGEN{ \mathbf{FinMon} \cap (-) },
    \]
where $\mathbf{FinMon}$ denotes the class of finite monoids;
\item the upper adjoint operator induced by $\mathsf{\pvV-Like}_\pvW$ is
    \[
        \CTSinPV{ \mathsf{\pvV-Like}_\pvW }{-} 
            \;=\; 
        ( - , \pvV ) \malcev \pvW;
    \]
\item the upper adjoint operator induced by $\mathsf{EPL}_\pvW$ is
    \[
        \CTSinPV{ \mathsf{EP}_\pvW }{-} 
            \;=\; 
        ( - ) \malcev \pvW;
    \]
\item the upper adjoint operator induced by $\mathsf{Ker}_\pvH$ is
    \[
        \CTSinPV{ \mathsf{Ker}_\pvH }{-} 
            \;=\; 
        ( - ) \malcev \pvH;
    \]
\item the upper adjoint operator induced by $\mathsf{EGen}$ is
    \[
        \CTSinPV{ \mathsf{EGen} }{-} 
            \;=\; 
        \mathbb{E}[-];
    \]
\item and the upper adjoint operator induced by $\mathsf{RegGen}$ is
    \[
        \CTSinPV{ \mathsf{RegGen} }{-} 
            \;=\; 
        \mathbb{R}[-].
    \]
\end{enumerate}
\end{example}

\begin{rmk}
Context specifiers are in some sense a generalization of the notion of "preimage classes" developed in \cite{CONGLATTICEPVAR}. 
They are also reminiscent of the notion of "divisor systems" developed in \cite{DIVISORSYSTEMS}.
The Galois connections featured in the results of these papers are---in the cases which overlap---the same as those featured here.
\end{rmk}

\begin{defn}
If $\scO$ is a context specifier, the \textbf{restriction} of a modulus $\Lambda$ to the context $\scO$ is the modulus
    \[
        \MDmodCTS{\Lambda}{\scO} 
            \;=\;
        \modrl{S}{ \bigcup_{U \in \scO(S)} \Lambda_U }.
    \]
The resulting map 
    \[
        \MDmodCTS{(-)}{\scO} : \MD \arw \MD
    \]
induces in turn a \textbf{context restriction} map
    \[
        \CFmodCTS{-}{\scO} : \CF \arw \CF 
    \]
whose action assigns $\cC \in \CF$ to the complex functor $\CFmodCTS{\cC}{\scO}$ given by
    \[
        \CFmodCTS{\cC}{\scO}(S)
            \;=\;
        \COMGEN{ \; \bigcup_{U \in \scO(S)} \cC(U) \; }{S}
    \]
at every finite semigroup $S$.
\end{defn}

\begin{thm} \label{cts:thm:moduli}
Let $\scO$ be a context specifier.
\begin{enumerate}
\item If $\Lambda$ is a modulus, then is $\MDmodCTS{\Lambda}{\scO}$ a modulus which refines $\Lambda$.
\item If $\Lambda_1$ and $\Lambda_2$ are moduli, then 
    \[
        \MDmodCTS{(\Lambda_1 \vee \Lambda_2)}{\scO}
            \;=\;
        (\MDmodCTS{\Lambda_1}{\scO}) \vee (\MDmodCTS{\Lambda_2}{\scO}).
    \]
\item The diagram
\begin{equation*}
            \begin{tikzcd}
                    \MD 
                    \ar[d, "\MDmodCTS{(-)}{\scO}"'] 
                    \ar[rr, "\points{-}"]
                &&
                    \PVAR 
                    \ar[d, "\CTSinPV{\scO}{-}"]
                \\
                    \MD 
                    \ar[rr, "\points{-}"']
                && 
                    \PVAR
            \end{tikzcd}
\end{equation*}
is commutative.
\end{enumerate}
\end{thm}

\begin{proof}
\hfill
\begin{enumerate}
\item
First, let $\varphi : S \rightarrow T$. 
If $X_S \in (\MDmodCTS{\Lambda}{\scO})_S$, then there exists some $U_S \in \scO(S)$ for which $X_S \in \Lambda_{U_S}$.
Now, there is some $U_T \in \scO(T)$ such that the image of $U_S$ under $\varphi$ is contained in $U_T$, and hence there exists some $X_T \in \Lambda_{U_T}$ such that $(X_S)\ext{\varphi} \subseteq X_T$.
Since $\Lambda_{U_T}$ is a subset of $(\MDmodCTS{\Lambda}{\scO})_T$, the first axiom is verified.

Next, let $\varphi : S \surj T$ be a regular epimorphism.
If $X_T \in (\MDmodCTS{\Lambda}{\scO})_T$, then there exists some $U_T \in \scO(T)$ for which $X_T \in \Lambda_{U_T}$.
Now, there is some $U_S \in \scO(S)$ such that $U_S$ maps onto $U_T$ under $\varphi$, and hence there exists some $X_S \in \Lambda_{U_S}$ such that $(X_S)\ext{\varphi} = X_T$.
Since $\Lambda_{U_S}$ is a subset of $(\MDmodCTS{\Lambda}{\scO})_S$, the second axiom is verified; and hence $\MDmodCTS{\Lambda}{\scO}$ is a modulus.

That $\MDmodCTS{\Lambda}{\scO}$ refines $\Lambda$ can be seen by considering for each finite semigroup $S$ the various inclusion maps $U \inj S$ as $U$ ranges over $\scO(S)$.
\item Straightforward.
\item A finite semigroup $S$ belongs to $\points{ \MDmodCTS{\Lambda}{\scO} }$ if and only if $\Lambda_U \subseteq \sing(U)$ for every subsemigroup $U \in \scO(S)$; that is, if and only if $\scO(S) \subseteq \points{\Lambda}$.
\end{enumerate}
\end{proof}

\begin{cor} \label{cts:cor:fptc}
Let $\scO$ be a context specifier.
The operator $\CFmodCTS{-}{\scO}$ is decreasing, preserves arbitrary joins, and satisfies FPTC with respect to $\CTSinPV{\scO}{-}$.
\end{cor}

\begin{proof}
The various claims follow from those of Theorem \ref{cts:thm:moduli} by way of Proposition \ref{mod:prop:construct}.
\end{proof}

\subsection{}
Of course, we are most intersted in the question of when Corollary \ref{cts:cor:fptc} can be strengthened to state that this pair satisfies PLTC.

\begin{example}[Varieties determined by subgroups]
The main result of \cite{PL-VARDETGRP-BEN} is equivalent to the statement that the monad completion of $\CFmodCTS{-}{\mathsf{Grp}}$ satisfies PLTC with respect to $\mathbb{G}[-]$; i.e., that
    \[
        \moncom{\CFmodCTS{\PV}{\mathsf{Grp}}} = \PL{ \mathbb{G}\pvV }
    \]
for every psuedovariety $\pvV$.
A more direct translation of Steinberg and van Gool's result is that $\moncom{\CFmodCTS{-}{\mathsf{Grp}}}$ satisfies the \textit{local} pointlike transfer condition with respect to $\mathbb{G}[-]$ on the interval $[ \pvtriv, \pvG]$---that is, that
    \[
        \moncom{\CFmodCTS{ \PL{\pvH} }{\mathsf{Grp}}} = \PL{ \mathbb{G}\pvH }
    \]
holds whenever $\pvH$ is a pseudovariety of groups.
The equivalence of these two statements follows from the fact that $\PV(G) = \PL{\pvV \cap \pvG}(G)$ whenever $G$ is a finite group.
More directly, the main result of \cite{PL-VARDETGRP-BEN} is that the modulus
    \[
        \MDmodCTS{ \mathsf{Ker}_\pvH }{ \mathsf{Grp} }
            \;=\;
        \modrl{S}{ \setst{ K_\pvH (G) }{ \text{$G$ is a subgroup of $S$} } }
    \]
is effective with respect to $\mathbb{G}\pvH$ for any pseudovariety of groups $\pvH$.
\end{example}

\begin{rmk}
A conjecture of Steinberg \cite[Conjecture~6.11]{PL-JOINS-BEN} (stated in the language of this paper) is that $\CFmodCTS{-}{ \mathsf{EP}_\pvW }$ satisfies PLTC with respect to $(-) \malcev \pvW$ whenever $\pvW$ is a pseudovariety of bands.
Note that since bands consist entirely of idempotents, in this case one has that
    \[
        \CFmodCTS{ \PV }{ \mathsf{EP}_\pvW }(S)
            \;=\;
        \COMGEN{ \; \bigcup \setst{ \, \PV \big( \langle X \rangle \big) }{ X \in \PL{\pvW}(S) \, } \; }{S}
    \]
for every pseudovariety $\pvV$ and every finite semigroup $S$. 

This conjecture is apparently inspired by his proof that the pseudovariety $\pvJ$ of $\GJ$-trivial semigroups has decidable pointlikes \cite[Theorem~6.9]{PL-JOINS-BEN}, wherein the decidability of $\pvJ$-pointlikes is derived from the equality
    \[
        \CFmodCTS{ \PL{\pvN} }{ \mathsf{EP}_\pvSL } = \PL{\pvN \malcev \pvSL}
    \]
by way of the fact that $\pvJ = \pvN \malcev \pvSL$.
\end{rmk}

\begin{prop}
Let $\scO$ be a context specifier.
If $S$ is a finite semigroup for which $S \in \scO(S)$ then the equality
    \[
        \PL{ \CTSinPV{\scO}{\pvV} }(S)       
            \;=\;   \PV(S)      
            \;=\;   \PL{\CTSsumPV{\scO}{\pvV}}(S)
    \]
holds for every psuedovariety $\pvV$.
\end{prop}

\begin{proof}
Let $\pvV$ be a pseudovariety.
The inequalities
    \[
        \PL{ \CTSinPV{\scO}{\pvV} }       
            \;\leq\;   \PV      
            \;\leq\;   \PL{\CTSsumPV{\scO}{\pvV}}
    \]
clearly hold.

Beginning with the converse of the left-hand inequality, the hypothesis that $S \in \scO(S)$ implies that
    \[
        \PV(S) 
            \;\subseteq\;
        \CFmodCTS{\PV}{\scO}(S)
            \;\subseteq\;
        \PL{ \CTSinPV{\scO}{\pvV} }(S)
    \]
by way of Corollary \ref{cts:cor:fptc}, from which it follows that $\PV(S) = \PL{ \CTSinPV{\scO}{\pvV} }(S)$.

To establish the converse of the right-hand inequality, let $\rho : S \rlm V$ be a relational morphism with $V \in \pvV$ which computes $\PV(S)$.
Then there exist $\scO$-subsemigroups $U_\rho \in \scO(\Gamma(S , \rho, V))$ and $U_V \in \scO(V)$ as in the diagram
    \begin{equation*}
    \begin{tikzcd}
    &&
        U_\rho
            \ar[dll, two heads, shift right=1]
            \ar[d, hook]
            \ar[rr]
    &&
        U_V
            \ar[d, hook]
        \\
        S
    &&
        \Gamma(S , \rho , V)
            \ar[ll, two heads]
            \ar[rr]
    &&
        V
    \end{tikzcd}
    \end{equation*}
yielding a relational morphism $\widetilde{\rho} : S \rlm U_V$ with graph $U_\rho$.
The fact that $U_V \in \CTSsumPV{\scO}{\pvV}$ yields the containment
    \[
        \PL{\CTSsumPV{\scO}{\pvV}}
            \;\subseteq\;
        \Nrv(S, \widetilde{\rho}, U_V)
            \;\subseteq\;
        \PV(S)
    \]
from which the proposition follows.
\end{proof}

\begin{cor}
Let $\pvV$ be a pseudovariety.
\begin{enumerate}
    \item If $S$ is generated by its idempotents, then
        \[
            \PL{ \mathbb{E}\pvV }(S) 
            \;=\; \PV(S) 
            \;=\; \PL{ \mathbf{eV} }(S),
        \]
    where $\mathbf{eV} = \CTSsumPV{\mathsf{EGen}}{\pvV}$.
    \item If $S$ is generated by its regular elements, then
        \[
            \PL{ \mathbb{R}\pvV }(S) 
            \;=\; \PV(S) 
            \;=\; \PL{ \mathbf{rV} }(S),
        \]
    where $\mathbf{rV} = \CTSsumPV{\mathsf{RegGen}}{\pvV}$.
    \item If $G$ is a group, then
        \[
            \PL{ \mathbb{G}\pvV }(G) 
                \;=\; \PV(G) 
                \;=\; \PL{ \pvV \cap \pvG }(G).
        \]
    \item If $M$ is a monoid, then
        \[
            \PL{ \mathbb{L}\pvV }(M) 
                \;=\; \PV(M) 
                \;=\; \PL{ \pvV \cap \mathbf{FinMon} }(M),
        \]
    where $\mathbf{FinMon}$ denotes the class of finite monoids.\footnotemark 
    \footnotetext{Here $\PL{ \pvV \cap \mathbf{FinMon} }$ is a notationally abusive way to denote the corresponding pointlike monad on the category of finite monoids.}
\end{enumerate}
\end{cor}

\begin{rmk}     \label{cts:rmk:GMC}
Many of the operators on $\PVAR$ considered here as examples of those induced by context specifiers are given in \cite[Section~2.4.3]{QTHEORY} as examples of continuous operators which do not satisfy the \textit{global Mal'cev condition} (hereafter referred to as \textit{GMC}), where an operator $\alpha$ satisfies GMC if 
    \[
        \pvW^\alpha 
            \;\leq\;
        ( \pvV^\alpha , \pvV ) \malcev \pvW
    \]
for all $\pvV , \pvW \in \PVAR$.
The \textit{raison d'{\^e}tre} for GMC is that continuous operators satisfying it are precisely those which arise from pseudovarieties of relational morphisms by way of \cite{QTHEORY}'s titular $\mathfrak{q}$ operator.
With this in mind, context specifiers provide a ``natural origin'' for a class of ubiquitous operators which do not neatly fit---relative to the similarly uniquitous class of operators which fit very nicely---into the framework established by \cite{QTHEORY}.
\end{rmk}

%

%% file: content/mainmatter/conclusion.tex
%

\section{Further remarks and open problems}      \label{section:conclusion}

This paper is intended to provide only the first steps in a much longer journey.
In this final section, we will describe a few potentially fruitful lines of inquiry enabled and suggested by the framework developed so far.

\subsection{On relational morphisms}
Our treatment of relational morphisms as \textit{objects} rather than as arrows is a departure from convention, but we believe it to be broadly advantageous. 
There are various classical results which might benefit from being presented from this point of view; the most notable of which being the derived category (or the derived semigroupoid) of a relational morphism.
Although we suspect that most of these "reframings" would yield at most aesthetic improvements, there is one notable exception to this suspicion: 
it would be very interesting to adopt this point of view in order to study various families of classes of relational morphisms (e.g. continuously closed classes, pseudovarieties, etc) in relation to the program outlined in \cite{QTHEORY}.

\subsection{Topology of semigroup complexes}   \label{problem:complexestopology}
Since semigroup complexes are precisely the semigroup objects in the category of finite abstract simplicial complexes, it would be interesting to investigate connections between their topological and algebraic properties.

As an initial illustration of what this would look like, consider the relation on the faces of a semigroup complex $(S, \cK)$ defined by
    \[
        X \sim Y 
            \;\; \Longleftrightarrow \;\;
        \text{$X$ and $Y$ are contained in precisely the same maximal faces.}
    \]
If this relation is an equivalence relation, then the induced relation on $S$ is a congruence whose equivalence classes are the maximal faces of $(S, \cK)$; moreover, $\cK$ is the nerve of the quotient map induced by said congruence.\footnotemark
\footnotetext{Carrying this line of thought a bit further yields a Galois connection between the lattice of congruences on $S$ and the lattice of $S$-complexes.}

\subsection{Nerve duality}
The nerve construction might lead one to define the \textit{co-nerve} of a relational morphism $\rho: S \rlm T$ to be the $T$-complex given by
    \[
        \mathsf{CoNrv}(S , \rho, T) 
            \;=\; 
        \COMGEN{ \, \setst{(s)\rho}{s \in S} \, }{T}
    \]
It would be very interesting to study ways to infer properties of the nerve from the co-nerve (and vice-versa).\footnotemark
\footnotetext{One such inference is apparent: if one of them is the singletons, then the maximal faces of the other are a partition.}

It is worth noting that this exact construction shows up in a wide variety of contexts (most of them non-algebraic).
In many situations, the simplicial complexes that I have termed the nerve and co-nerve are called the \textit{Dowker complexes} of a relation---this is due to their consideration in \cite{DOWKERHOMREL}, wherein they are shown to have isomorphic homology and cohomology groups; and, moreover, their geometric realizations are shown to have the same homotopy type.\footnotemark
\footnotetext{A nice version of the proof may be found on the nLab at \cite{nlab:dowker}.}
A particular application of this is to the nerve and Vietoris complexes associated to a covering of a space---these are Dowker complexes, and their having isomorphic homology and cohomology is a particular case of the more general statement regarding Dowker complexes.
Hence this is also true for our nerve and co-nerve. 
I don't see any immediate applications---but sufficient progress on \ref{problem:complexestopology} might provide an opportunity for one.

\subsection{Lattice theoretic aspects of complex functors}
It will be necessary to investigate the lattice of complex functors from a lattice-theoretic point of view. 
A reasonable starting point would be to characterize various "special" classes of complex functors (e.g. those which are compact, co-compact, join/meet irreducible, and so on). 
In particular, relationships between complex functors satisfying these properties and the lattice-theoretic properties of their respective fixed point pseudovarieties should be considered.
Additionally, it will be important to study operators on $\CF$ with an eye towards applications to fixed point and pointlike transfer conditions. 
Material found in \cite{QTHEORY} is likely to provide inspiration for fruitful lines of inquiry in this realm.

\subsection{False pointlikes}
Remark \ref{existenceofimposters} shows that the Galois connection between $\CM$ and $\PVAR$ is not an equivalence by giving an example---due originally to Almeida \textit{et~al.} in \cite{PL-RJ-ALM}---of a complex monad which is \textit{not} one of the pointlike monads.
For ease of reference---and for the sake of an evocative name---let us call a complex monad which fixes $\pvV$ but which is strictly below $\PV$ an \textit{imposter} with respect to $\pvV$.
The fact that there exists a case (that of $\pvJ$) for which at least one imposter exists is the extent of our knowledge about them.
Otherwise, this phenomenon is mysterious. 
Are there pseudovarieties for which no imposters exist? 
When they do exist, how many of them can there be? Are there pseudovarieties with finitely many, countably many, continuum many imposters? 
What can information on imposters tell you about the properties of a pseudovariety (and vice-versa)? 
Also, are there classes of semigroups within which ``locally'' no imposters exist, i.e., on which the actions of all complex monads fixing a given pseudovariety coincide?

\subsection{Pointlikes for semidirect products}
A vitally important notion in finite semigroup theory which is conspicuously absent here is that of \textit{semidirect products}.
A natural program would be to seek operators on $\CF$ whose pairing with operators of the form $\pvV \ast (-)$ and $(-) \ast \pvV$ satisfies the fixed point transfer condition.

One approach to this would involve developing an analogue of context specifiers which is related to the notion of $\pvV$-stabilizer pairs. 
Recall that a \textit{$\pvV$-stabilizer pair} of a finite semigroup $S$ is a pair $(x, U)$ where $x \in S$ and $U$ is a subsemigroup of $S$ such that for any relational morphism $\rho : S \rlm V$ with $V \in \pvV$ there exists an element $v \in V$ such that $x \in (v)\rho^\inv$ and $U$ is a subsemigroup of the inverse image of the right stabilizer of $v$.
It is shown in \cite[Lemma~5.3]{APRELMOR} that if $\pvW$ is a Tilson-local pseudovariety then
    \[
        \setst{ x \cdot A}{A \in \PL{\pvW}(U)}
            \;\subseteq\;
        \PL{\pvW \ast \pvV}(S)
    \]
whenever $(x, U)$ is a $\pvV$-stabilizer pair of $S$.
This suggests a generalization involving evaluating some modulus on subsemigroups parameterized by particular semigroup elements, then taking products of the resulting pairs.

Another approach would be to generalize the entire framework of this paper to the context of finite semigroupoids.
The idea here would be to build on the results of \cite{PL-DELAYTHM} and \cite{PL-JOINS-BEN}, wherein the notion of pointlike sets for finite semigroupoids (and finite categories) are developed in pursuit of pointlike results for semidirect products of pseudovarieties.
It was shown in \cite{PL-DELAYTHM} that pointlikes retain the analogous properties in these contexts to those which were generalized in the context of finite semigroups to yield the framework of this paper---not only are they functors that preserve embeddings and quotient functors; but they are in fact \textit{monads} whose multiplication transformation is the union map.
Hence it is very likely that the whole story plays out in the context of semigroupoids in precisely the same manner; the definitions basically write themselves.
Once the analogous notions are established, the evident goal would be to find an adjoint triple between the lattices of complex functors for finite semigroupoids and finite semigroups which "mirrors" the adjoint triple 
    \[
        (\mathbf{g} \; \dashv \;  (-) \cap \FSGP \; \dashv \; \ell) : \PVAR \arw \mathbb{SDPV},
    \]
where $\mathbb{SDPV}$ denotes the lattice of pseudovarieties of finite semigroupoids and $(-) \cap \FSGP : \mathbb{SDPV} \rightarrow \PVAR$ is a notationally abusive way of denoting the map assigning a semigroupoid pseudovariety to the semigroup pseudovariety whose members are the endomorphism semigroups of its one-object members.

%